\documentclass[a4paper,12pt,reqno]{amsart}
\usepackage[T1]{fontenc}
\usepackage[utf8]{inputenc}
\usepackage[margin=1in]{geometry}
\usepackage[libertinus,vvarbb]{newtx}
\usepackage{enumitem}
\usepackage{graphicx,color}
\usepackage{xparse}
\usepackage[bookmarksdepth=2]{hyperref}

\newcommand{\ignore}[1]{}

\numberwithin{equation}{section}

\newtheorem{theorem}{Theorem}

\newtheorem{corollary}[theorem]{Corollary}
\newtheorem{proposition}[theorem]{Proposition}

\theoremstyle{definition}

\newtheorem{problem}{}[section]

\DeclareMathOperator{\id}{Id}

\DeclareMathOperator{\artanh}{ar\,tanh}

\newcommand{\ind}{\mathbb{1}}

\newcommand{\ex}{\mathbb{E}}
\newcommand{\C}{\mathbb{C}}
\newcommand{\R}{\mathbb{R}}
\newcommand{\Z}{\mathbb{Z}}

\newcommand{\fourier}{\mathscr{F}}
\newcommand{\ph}{\varphi}
\newcommand{\eps}{\varepsilon}

\newcommand{\pvint}{\operatorname{pv}\!\int}
\newcommand{\dir}{{\mathrm{D}}}
\newcommand{\neu}{{\mathrm{N}}}

\renewcommand{\le}{\leqslant}
\renewcommand{\ge}{\geqslant}

\NewDocumentCommand{\formula}{ssom}{%
 \IfBooleanTF{#1}{%
  \IfBooleanTF{#2}{%
   \IfValueTF{#3}%
    {\begin{align}\label{#3}\begin{gathered}#4\end{gathered}\end{align}}%
    {\begin{gather}#4\end{gather}}%
  }{%
   \IfValueTF{#3}%
    {\begin{align}\label{#3}\begin{aligned}#4\end{aligned}\end{align}}%
    {\begin{gather*}#4\end{gather*}}%
  }%
 }{%
  \IfValueTF{#3}%
   {\begin{align}\label{#3}#4\end{align}}%
   {\begin{align*}#4\end{align*}}%
 }%
}

\newcommand{\doi}[1]{\href{https://doi.org/#1}{\textsf{\scriptsize DOI:#1}}}
\newcommand{\arxiv}[1]{\href{https://arxiv.org/abs/#1}{\textsf{\scriptsize arXiv:#1}}}

\begin{document}

\title[Harmonic extension technique]{Harmonic extension technique: \\ probabilistic and analytic perspectives}
\author{Mateusz Kwaśnicki}
\address{Mateusz Kwaśnicki \\ Department of Pure Mathematics \\ Wrocław University of Science and Technology \\ ul. Wybrzeże Wyspiańskiego 27 \\ 50-370 Wrocław, Poland}
\email{\href{mailto:mateusz.kwasnicki@pwr.edu.pl}{mateusz.kwasnicki@pwr.edu.pl}}

\begin{abstract}
Consider a path of the reflected Brownian motion in the half-plane $\{y \ge 0\}$, and erase its part contained in the interior $\{y > 0\}$. What is left is, in an appropriate sense, a path of a jump-type stochastic process on the line $\{y = 0\}$ --- the boundary trace of the reflected Brownian motion. It is well known that this process is in fact the $1$-stable Lévy process, also known as the Cauchy process.

The PDE interpretation of the above fact is the following. Consider a bounded harmonic function $u$ in the half-plane $\{y > 0\}$, with sufficiently smooth boundary values~$f$. Let $g$ denote the normal derivative of $u$ at the boundary. The mapping $f \mapsto g$ is known as the Dirichlet-to-Neumann operator, and it is again well known that this operator coincides with the square root of the 1-D Laplace operator $-\Delta$. Thus, with a suitable normalisation of boundary local time, the negative of the Dirichlet-to-Neumann operator is the generator of the boundary trace process.

Molchanov and Ostrovskii proved that isotropic stable Lévy processes are boundary traces of appropriate diffusions in half-spaces. Caffarelli and Silvestre gave a PDE counterpart of this result: the fractional Laplace operator is the Dirichlet-to-Neumann operator for an appropriate second-order elliptic equation in the half-space. Again, after matching the normalisations, the negative of the Dirichlet-to-Neumann operator is the generator of the boundary trace process.

During my talk I will discuss boundary trace processes and Dirichlet-to-Neumann operators in a more general context. My main goal will be to explain the connections between probabilistic and analytical results. Along the way, I will introduce the necessary machinery: Brownian local times and additive functionals, Krein's spectral theory of strings, and Fourier transform methods.
\end{abstract}

\maketitle
\thispagestyle{empty}

%
%

These are the lecture notes for a series of talks given at the \emph{3rd Korean Croatian Summer Probability Camp} in Šibenik, Croatia, on June 26--28, 2024. I sincerely thank the organisers, \textsc{Panki Kim, Nikola Sandrić, Ana Perišić,} and \textsc{Vanja Wagner,} for letting me be a part of this lovely event.

I aim to keep the exposition conversational, focusing on core ideas rather than technical formalities. As a result, some statements are deliberately informal. Unless stated otherwise, all functions of the spatial variables $x, y, z$ are real-valued, and boundary data are continuous and compactly supported. Occasionally, integrals appearing in identities are implicitly assumed to be finite. Finally, generic positive constants like $c$ or $c_{d, \alpha}$ may vary from line to line, or even within a single equation.

If you find a mistake, have a comment, or simply find these notes useful, please do send me a message at \href{mailto:mateusz.kwasnicki@pwr.edu.pl}{\textsf{mateusz.kwasnicki@pwr.edu.pl}}\textsf{.}

\newpage
\section{Cauchy process and the Dirichlet-to-Neumann operator}

\subsection{Dirichlet-to-Neumann map}

\begin{figure}
\noindent
\hfill
\begin{minipage}{0.5\textwidth}
\centering
\includegraphics[page=1,scale=0.75]{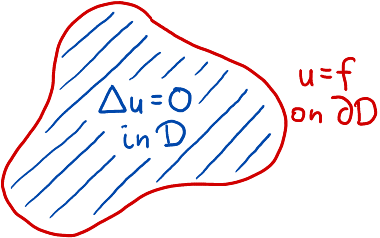}
\caption{The Dirichlet problem.}
\label{fig:dirichlet}
\end{minipage}%
\begin{minipage}{0.5\textwidth}
\centering
\includegraphics[page=2,scale=0.75]{pictures.pdf}
\caption{Normal derivative at the boundary.}
\label{fig:dtn}
\end{minipage}
\hfill
\end{figure}

A function $u$ on a Euclidean domain $D \subseteq \R^d$ is harmonic if it is a solution of the Laplace equation $\Delta u = 0$, where, as usual,
\[
 \Delta u = \sum_{j = 1}^d \frac{\partial^2 u}{\partial x_j^2}
\]
is the Laplace operator. If $D$ is nice enough (for example, a bounded smooth domain), then for every continuous function $f$ on $\partial D$ there is a unique solution of the Dirichlet problem (see Figure~\ref{fig:dirichlet})
\[
 \begin{cases}
  \Delta u = 0 & \text{in $D$,} \\
  u = f & \text{on $\partial D$.}
 \end{cases}
\]
Suppose that $f$ is smooth enough (for example, $C^\infty$), and let $g$ denote the (outward) normal derivative of $u$ --- the corresponding Neumann boundary condition. Then $g$ is given in terms of $f$ by a certain linear operator $K$:
\[
 g = K f ,
\]
known as the Dirichlet-to-Neumann operator.

This operator, in a slightly more general form, arises in Calderón's inverse conductivity problem: knowing only the boundary measurements, can one reconstruct the conductivity inside $D$? In this case the Laplace operator $\Delta$ is replaced by $\nabla \cdot (\gamma \nabla u)$, the normal derivative is replaced by the conormal derivative $\gamma \partial_\nu u$, and the question asks whether knowledge of this modified $K$ determines $\gamma$.

If $P_D(x, x')$ denotes the Poisson kernel of $D$ and $\sigma$ is the surface measure on $\partial D$, then
\[
 u(x) = \int_{\partial D} f(x') P_D(x, x') \sigma(dx') .
\]
It is thus reasonable to expect that the kernel of $K$ is given by the outward normal derivative at the boundary of the Poisson kernel; say,
\[
 K f(x) = -\int_{\partial D} f(x') J(x, x') \sigma(dx')
\]
with
\[
 J(x, x') = -\frac{\partial P_D}{\partial \nu_x}(x, x') = \frac{\partial^2 G_D}{\partial \nu_x \partial \nu_{x'}}(x, x') ,
\]
where $\nu$ denotes the outward normal vector and $G_D(x, x')$ is the Green function for $-\Delta$ (see Figure~\ref{fig:dtn}; forget about it if you are not already familiar with the Green function and the Poisson kernel, this paragraph is only meant to put things in a broader context). This is plainly wrong: $K$ applied to a constant function is zero, so we would rather have something like
\[
 K f(x) = \int_{\partial D} (f(x) - f(x')) J(x, x') \sigma(dx') ,
\]
with the extra $f(x)$ term originating from the singularity of $P_D(x, x')$ as $x$ approaches $x'$. The above expression turns out to be correct, as long as the integral is understood in the Cauchy principal value sense, but we will not pursue this direction.

Before we proceed, let us mention another description of the operator $K$. With the notation introduced above, by Green's first identity,
\[
 \int_D |\nabla u(x)|^2 dx = \int_{\partial D} u(x) \, \frac{\partial u}{\partial \nu}(x) \sigma(dx) = \int_{\partial D} f(x) g(x) \sigma(dx) .
\]
Dirichlet's principle states that $u$ is the unique minimiser of the integral on the left-hand side of the above equation among sufficiently regular functions on $\overline D$ with boundary values $f$. Thus,
\[
 \int_{\partial D} f(x) K f(x) \sigma(dx) = \min \biggl\{ \int_D |\nabla v(x)|^2 dx \biggr\} ,
\]
where the minimum is taken over the class of admissible functions $v$ just described, and it is attained when $v = u$ is the harmonic function with boundary values $f$.

\subsection{Square root of the Laplace operator}

\begin{figure}
\noindent
\hfill
\begin{minipage}{0.5\textwidth}
\centering
\includegraphics[page=3,scale=0.75]{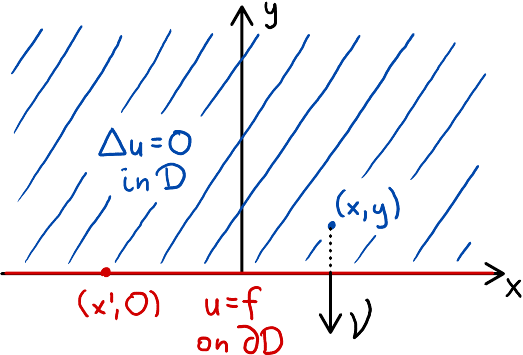}
\caption{The half-plane.}
\label{fig:hp}
\end{minipage}%
\begin{minipage}{0.5\textwidth}
\centering
\includegraphics[page=4,scale=0.75]{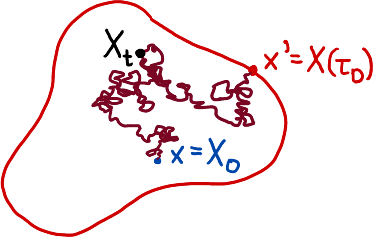}
\caption{Brownian motion and the first exit time.}
\label{fig:kakutani}
\end{minipage}
\hfill
\end{figure}

Throughout these lecture notes we will study the Dirichlet-to-Neumann operator for the $(d + 1)$-dimensional half-space $D = \R^d \times (0, \infty)$, the boundary of which is naturally identified with $\R^d$. This domain is unbounded, and so in order to have a unique solution of the Dirichlet problem, we need to assume, for example, that $f$ is compactly supported and $u$ is bounded. From now on we use the symbol $x \in \R^d$ for the `horizontal' coordinate and $y \in (0, \infty)$ for the vertical one (see Figure~\ref{fig:hp}).

The Poisson kernel is given by
\[
 P_D((x, y), x') = P_y(x - x') = c_d \, \frac{y}{(|x - x'|^2 + y^2)^{(d + 1) / 2}}
\]
(to avoid nested brackets, from now on we switch the notation to $P_y(x - x')$ when speaking about the Poisson kernel for the half-space), the harmonic extension becomes a convolution with $P_y$:
\[
 u(x, y) = \int_{\R^d} f(x') P_y(x - x') dx' = f * P_y(x) ,
\]
and it is more or less straightforward to see that
\[
 K f(x) = -\frac{\partial u}{\partial y}(x, 0) = c_d \pvint_{\R^d} \frac{f(x) - f(x')}{|x - x'|^{d + 1}} \, dx' .
\]
This operator is sometimes called half-Laplacian and denoted by $\sqrt{-\Delta}$. In the sense of spectral theory, it is indeed the square root of $-\Delta$, the $d$-dimensional Laplace operator (with a minus sign). A somewhat informal way to prove this is as follows: since the $(d + 1)$-dimensional Laplace operator of $u$ is zero, we have
\[
 K^2 f(x) = \frac{\partial^2 u}{\partial y^2}(x, 0) = -\sum_{j = 1}^d \frac{\partial^2 u}{\partial x_j^2}(x, 0) = -\Delta f(x) ;
\]
and by Dirichlet's principle, $K$ is non-negative definite. A rigorous argument involves the $d$-dimensional Fourier transform:
\[
 \fourier f(\xi) = \int_{\R^d} e^{-i \xi x} f(x) dx
\]
(here $\xi x$ is the usual dot product). If we denote by $\fourier u(\xi, y)$ the Fourier transform of $u(\cdot, y)$, then the Laplace equation is turned into an ODE
\[
 -|\xi|^2 \fourier u(\xi, y) + \frac{\partial^2 \fourier u}{\partial y^2}(\xi, y) = 0 ,
\]
which, given that $u$ is assumed to be bounded, can be shown to have a unique solution, given by
\[
 \fourier u(\xi, y) = e^{-|\xi| y} \fourier u(\xi, 0) = e^{-|\xi| y} \fourier f(\xi) .
\]
Note that this actually shows that $e^{-|\xi| y}$ is the Fourier transform of the Poisson kernel $P_y(\cdot)$. It follows that
\[
 \fourier K f(\xi) = -\frac{\partial \fourier u}{\partial y}(\xi, 0) = |\xi| \fourier u(\xi, 0) = |\xi| \fourier f(\xi) .
\]
In particular, $\fourier K^2 f(\xi) = |\xi|^2 \fourier f(\xi) = -\fourier \Delta f(\xi)$.

\subsection{Boundary trace process}

Up to a factor $\tfrac12$, the Laplace operator is the generator of the Brownian motion. The $d$-dimensional Brownian motion $X_t$ is just a vector of independent Wiener processes, and the previous sentence means that
\[
 \tfrac12 \Delta u(x) = \lim_{t \to 0^+} \frac{\ex [u(X_t) \vert X_0 = x] - u(x)}{t}
\]
for sufficiently regular functions $u$. The connection is, however, much deeper. For example, Kakutani observed that the Poisson kernel $P_D(x, \cdot)$ is the density function (with respect to the surface measure) of the distribution of $X_t$ at the first exit time from $D$. Thus, if we denote this exit time by $\tau_D$, the solution of the Dirichlet problem is given by (see Figure~\ref{fig:kakutani})
\[
 u(x) = \ex [f(X(\tau_D)) \vert X_0 = x] .
\]
Note that we occasionally write $X(t)$ instead of $X_t$ to avoid nested subscripts.

Let us consider a $(d + 1)$-dimensional Brownian motion $(X_t, Y_t)$ in the half-space $D = \R^d \times [0, \infty)$, reflected at the boundary, and started at $(X_0, Y_0) = (0, 0)$ unless stated otherwise. Formally, this can be constructed from an unrestricted $(d + 1)$-dimensional Brownian motion $(X_t, W_t)$ by setting $Y_t = |W_t|$.

\begin{figure}
\noindent
\hfill
\begin{minipage}{0.5\textwidth}
\centering
\includegraphics[page=5,scale=0.75]{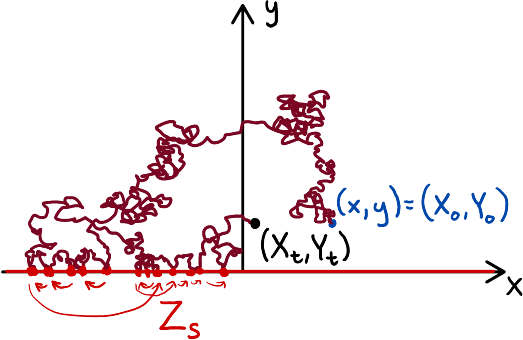}
\caption{Boundary trace.}
\label{fig:trace}
\end{minipage}%
\begin{minipage}{0.5\textwidth}
\centering
\includegraphics[page=6,scale=0.75]{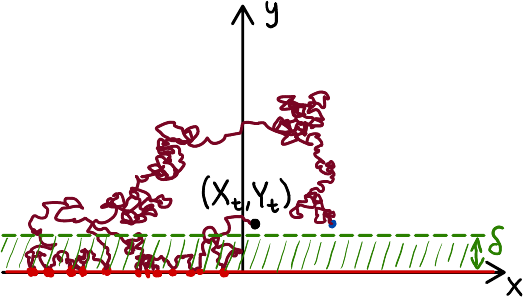}
\caption{Limiting procedure.}
\label{fig:limit}
\end{minipage}
\hfill
\end{figure}

What is the trace left by $(X_t, Y_t)$ on the boundary? Well, it is just a collection of dots --- an uncountable one, with no isolated points, almost surely (see Figure~\ref{fig:trace}). But is there a better way to view this trace? It turns out that the answer is yes: we can identify the boundary trace with a certain stochastic process on $\R^d$. This, however, requires some care, as $(X_t, Y_t)$ spends zero time at the boundary.

What we are going to do is the following: we take out all fragments of the paths of $(X_t, Y_t)$ which lie outside $\R^d \times [0, \delta]$ (without gluing the pieces together: the resulting process has jumps), we scale the time appropriately, and we hope that in the limit as $\delta \to 0^+$ we will get a nontrivial process $Z_t$ on the boundary (see Figure~\ref{fig:limit}). A rigorous construction involves what is called the local time of $Y_t$ at $0$.

\begin{figure}
\noindent
\hfill
\begin{minipage}{0.5\textwidth}
\centering
\includegraphics[page=7,scale=0.75]{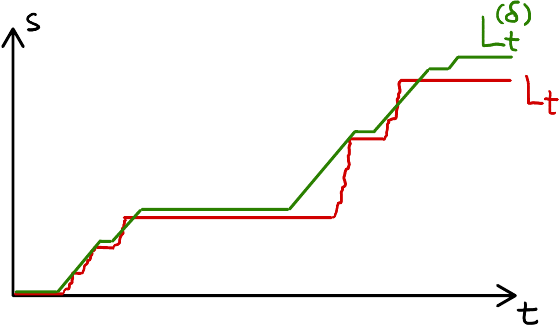}
\caption{Local time.}
\label{fig:local}
\end{minipage}%
\begin{minipage}{0.5\textwidth}
\centering
\includegraphics[page=8,scale=0.75]{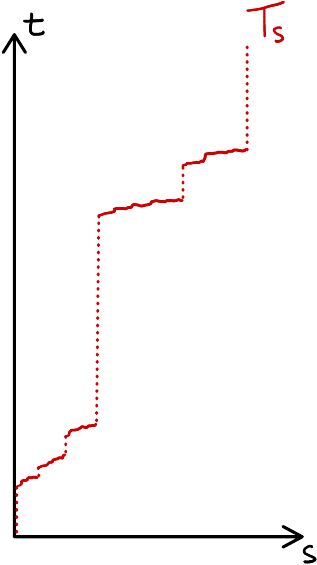}
\caption{Inverse local time.}
\label{fig:inverse}
\end{minipage}
\hfill
\end{figure}

For a given $\delta > 0$ and $t \ge 0$ let $L_t^{(\delta)}$ denote the amount of time spent by the process $(X_s, Y_s)$ up to time $t$ in the $\delta$-neighbourhood of the boundary:
\[
 L_t^{(\delta)} = \int_0^t \ind_{[0, \delta]}(Y_s) ds .
\]
Of course, $L_t^{(\delta)}$ converges to zero as $\delta \to 0^+$. However, it can be proved that $(2 \delta)^{-1} L_t^{(\delta)}$ converges with probability one to a nontrivial continuous process $L_t$, called the local time of $Y_t$ at $0$. Clearly, $L_t$ is non-decreasing, and it is constant on every interval of time where $Y_t > 0$. The nontrivial part is that $L_t$ is finite and nonzero; more precisely, almost surely every $t \ge 0$ such that $Y_t = 0$ is a point of increase of $L_t$. So, in other words, $dL_t$ is a continuous measure whose support is the set of times $t$ such that $Y_t = 0$ (see Figure~\ref{fig:local}).

The local time is not invertible, but it admits a generalised right-continuous inverse
\[
 T_s = \max\{t \ge 0 : L_t \le s\} ,
\]
called simply the inverse local time of $Y_t$ at $0$, and commonly denoted by $L_s^{-1}$ (see Figure~\ref{fig:inverse}). Then $L(T_s) = s$ for every $s \ge 0$, but we only have $T(L_t) \ge t$ for $t \ge 0$. Furthermore, with probability one the range of the process $T_s$ (augmented by the left limits $T_{s-}$, to be precise) coincides with the set of $t$ such that $Y_t = 0$. In particular, $Y(T_s) = 0$.

This tells us that the process
\[
 Z_s = X(T_s)
\]
(again: augmented by left limits $Z_{s-}$) is precisely what we are after: it is the trace left by $(X_t, Y_t)$ on the boundary.

\subsection{Cauchy process}

We already know how we can construct the boundary trace process $Z_s$, but what kind of stochastic process is it? An answer to this question was found by Spitzer in 1958~\cite{spitzer}.

It turns out that $T_s$ is a Lévy process, that is, a càdlàg (right-continuous with left limits) process with stationary and independent increments. Intuitively, this is a consequence of the strong Markov property of $Y_t$: the magnitudes of jumps of $T_s$ are the lengths of excursions of $Y_t$ away from $0$. However, a rigorous proof of this property is more complicated. Since $T_s$ is non-decreasing, by definition it is a subordinator.

Since $X_t$ and $T_s$ are clearly independent processes, $X_t$ is the Brownian motion, and $T_s$ is a subordinator, we find that $Z_s$ is a subordinate Brownian motion. Furthermore, it is known that $T_s$ is the $\tfrac12$-stable subordinator, and so it follows that $Z_s$ is the isotropic $1$-stable Lévy process in $\R^d$, known under the name Cauchy process.

The Cauchy process is perhaps the best understood multidimensional jump-type process. The random variables $Z_s$ are known to have Cauchy distributions (which explains the name of the process), and hence the density function of $Z_s$ coincides with the Poisson kernel for the half-space. The generator of $Z_s$ is $-\sqrt{-\Delta}=-K$, the negative of the Dirichlet-to-Neumann operator. The Lévy measure of $Z_s$, which describes the intensity of jumps of $Z_s$ of a given magnitude, has a density function $c_d |z|^{-d - 1}$, identical to the non-local kernel of the operator $K$. Are these properties merely a coincidence, or is there a deeper theory behind them?

In the next section we study a particular modification of what has been discussed above, which already clarifies the situation. A general answer to the above questions is the main result discussed in these notes. But let us make it clear straight away: the first property (equality of the density function of $Z_s$ and the Poisson kernel) is a serendipity, but the other two manifest certain general principles.

\subsection*{Exercises}

\begin{problem}
Let $K$ be the Dirichlet-to-Neumann operator for the half-space $\R^d \times (0, \infty)$. Prove that $K$ is homogeneous with degree $1$: if $\lambda > 0$ and $f_\lambda(x) = f(\lambda x)$, then $K f_\lambda(x) = \lambda K f(\lambda x)$.
\end{problem}

\begin{problem}
\begin{enumerate}[label={\textup{(\alph*)}}]
\item
Verify directly that the Poisson kernel for the half-space $D = \R^d \times (0, \infty)$:
\[
 P_y(x - x') = c_d \, \frac{y}{(|x - x'|^2 + y^2)^{(d + 1) / 2}} \, ,
\]
is indeed a harmonic function of $(x, y)$.
\item
Prove that there exists a constant $c_d$ such that
\[
 \int_{\R^d} P_y(x - x') dx' = 1 .
\]
\item
With this choice of $c_d$, consider a compactly supported continuous function $f$, and denote
\[
 u(x, y) = \int_{\R^d} P_y(x - x') f(x') dx' .
\]
Prove that $u$ is indeed the harmonic extension of $f$: a bounded solution to the Dirichlet problem in $D$ with boundary values $f$.
\item
Bonus problem: Find a closed-form expression for $c_d$.
\end{enumerate}
\end{problem}

\begin{problem}
Let $f$ be a compactly supported smooth function on $\R^d$. Prove rigorously that the Dirichlet-to-Neumann operator for the half-space is indeed given by
\[
 \begin{aligned}
  K f(x) & = \frac{c_d}{2} \int_{\R^d} \frac{2 f(x) - f(x + z) - f(x - z)}{|z|^{d + 1}} \, dz \\
  & = c_d \pvint_{\R^d} \frac{f(x) - f(x')}{|x - x'|^{d + 1}} \, dx' .
 \end{aligned}
\]
\end{problem}

\begin{problem}
Let $f$ be a compactly supported continuous function on $\R^d$, and let $u$ be a bounded harmonic function on $\R^d \times (0, \infty)$ with boundary values $f$. Prove rigorously that $\fourier u(\xi, y) = e^{-|\xi| y} \fourier f(\xi)$ for every $\xi \in \R^d$ and $y > 0$.
\end{problem}

\begin{problem}
Let $Y_t$ be the reflected Brownian motion in $[0, \infty)$, and let $L_t$ be the local time of $Y_t$ at zero, defined as the limit of $(2 \delta)^{-1} L_t^{(\delta)}$. Assume that $Y_0 = 0$.
\begin{enumerate}[label={\textup{(\alph*)}}]
\item
Using self-similarity of $Y_t$, prove that $L_t$ is self-similar with Hurst exponent~$\tfrac12$: the process $L_{\lambda t}$ has the same law as the process $\lambda^{1/2} L_t$.
\item
Observe that the inverse local time $T_s$ is thus self-similar with Hurst exponent $2$, and therefore it is the $\tfrac12$-stable subordinator.
\item
Conclude that the boundary trace process $Z_s = X(T_s)$ is self-similar with Hurst exponent $1$, and therefore it is the isotropic $1$-stable Lévy process (the Cauchy process).
\end{enumerate}
\end{problem}

%
%

\newpage
\section{Fractional Laplacian and isotropic stable Lévy processes}

\subsection{Fractional Laplace operator}

There is little special about the half-Laplacian among other fractional powers of the Laplace operator. These can be defined in terms of the Fourier transform:
\[
 \fourier (-\Delta)^{\alpha/2} f(\xi) = |\xi|^\alpha \fourier f(\xi) ,
\]
or as a principal value integral:
\[
 (-\Delta)^{\alpha/2} f(x) = c_{d, \alpha} \pvint_{\R^d} \frac{f(x) - f(x')}{|x - x'|^{d + \alpha}} \, dx' .
\]
The former definition works for arbitrary $\alpha > 0$; in fact, complex powers can also be defined. The latter one is more common, and it requires $\alpha \in (0, 2)$ (although a similar expression involving a hyper-singular integral can be given for more general $\alpha$). From now on, $\alpha \in (0, 2)$ is our standing assumption.

Can one obtain $(-\Delta)^{\alpha/2}$ as a Dirichlet-to-Neumann operator? Recall that the jump kernel associated with $K$ was essentially the normal derivative of the Poisson kernel. Thus, a naive way to go would be to replace the $(d + 1)$-dimensional Laplace operator by something else, so that the analogue of the Poisson kernel for the half-space $D = \R^d \times (0, \infty)$ would be, say,
\[
 P_y(x - x') = c_{d, \alpha} \, \frac{y}{(|x - x'|^2 + y^2)^{(d + \alpha) / 2}} \, .
\]
This fails for obvious reasons: as $y \to 0^+$, the integral of $P_y(\cdot)$ converges to $0$ (if $\alpha < 1$) or to $\infty$ (if $\alpha > 1$), while we require it to converge to $1$. To remedy this, we can try playing with the exponents, and a simple scaling argument brings us to
\[
 P_y(x - x') = c_{d, \alpha} \, \frac{y^\alpha}{(|x - x'|^2 + y^2)^{(d + \alpha) / 2}} \, .
\]
If $c_{d, \alpha}$ is chosen appropriately, then $P_y(\cdot)$ turns out to have integral $1$ for every $y > 0$, and so indeed $(-\Delta)^{\alpha/2}$ can be identified with an operator reminiscent of the Dirichlet-to-Neumann map: it is more or less clear that the fractional Laplace operator $(-\Delta)^{\alpha/2}$ can be expressed as
\[
 \begin{aligned}
  K f(x) & = -\lim_{y \to 0^+} \frac{u(x, y) - u(x, 0)}{y^\alpha} \\
  & = -\lim_{y \to 0^+} \frac{1}{\alpha y^{\alpha - 1}} \, \frac{\partial u}{\partial y}(x, y) ,
 \end{aligned}
\]
where
\[
 u(x, y) = \int_{\R^d} f(x') P_y(x - x') dx'
\]
is an extension of $f$ to the half-space $\R^d \times (0, \infty)$. (This last statement is slightly wrong: $K$ is not exactly equal to $(-\Delta)^{\alpha/2}$. Can you spot the error?)

Perhaps the most surprising aspect of the above construction is that it actually works, in the sense that it is linked to a certain partial differential equation in the half-space. One can directly verify that the Poisson kernel $P_y(x - x')$, and therefore also the extension $u$, satisfies the following modified Laplace equation:
\[
 \Delta_{x, y} u(x, y) + \frac{1 - \alpha}{y} \, \frac{\partial u}{\partial y}(x, y) = 0 ,
\]
which can be written as an equation in divergence form:
\[
 \nabla_{x, y} \cdot (y^{1 - \alpha} \nabla_{x, y} u)(x, y) = 0 .
\]
That is, $u(x, y)$ is `harmonic' with respect to the operator
\[
 L = \Delta_{x, y} + \frac{1 - \alpha}{y} \, \frac{\partial}{\partial y} \, ,
\]
or, equivalently, with respect to
\[
 y^{1 - \alpha} L = \nabla_{x, y} \cdot (y^{1 - \alpha} \nabla_{x, y}) .
\]
While the above construction may seem very natural, it became widely used following the seminal work by Caffarelli and Silvestre from 2007~\cite{caffarelli-silvestre}, and it is known as the Caffarelli--Silvestre extension technique.

We leave the details of the above calculations aside, and instead we discuss an alternative solution that involves the Fourier transform. It turns the PDE $L u = 0$ into an ODE
\[
 -|\xi|^2 \fourier u(\xi, y) + \frac{\partial^2 \fourier u}{\partial y^2}(\xi, y) + \frac{1 - \alpha}{y} \, \frac{\partial \fourier u}{\partial y}(\xi, y) = 0 .
\]
For $\xi \ne 0$ this is a variant of the Bessel differential equation; boundedness and the prescribed value at $0$ select a unique solution. It follows that
\[
 \fourier u(\xi, y) = c_\alpha (y |\xi|)^{\alpha/2} K_{\alpha/2}(y |\xi|) \fourier f(\xi) ,
\]
where $K_{\alpha/2}$ is the modified Bessel function of the second kind, and one can verify directly that
\[
 \fourier K f(\xi) = -\lim_{y \to 0^+} \frac{\fourier u(\xi, y) - \fourier u(\xi, 0)}{y^\alpha} = c_\alpha |\xi|^\alpha \fourier f(\xi)
\]
for appropriate constants $c_\alpha$ (different from the one in the previous equation; we do not pay much attention to particular values of constants in these notes).

The Caffarelli--Silvestre extension technique allows one to convert certain non-local problems in $\R^d$ to local problems in $\R^d \times (0, \infty)$. Many results of the theory of PDEs are available only for local problems, either because they generally fail in the non-local setting, or because the methods to tackle their non-local variants have not been developed yet. In these cases the Caffarelli--Silvestre extension technique is the right way to go.

The PDE discussed above appeared already in an article by Muckenhoupt and Stein from 1965~\cite{muckenhoupt-stein}, but in a much different context. Quite surprisingly, however, a probabilistic counterpart of the above construction goes back to the paper by Molchanov and Ostrovskii from 1969~\cite{molchanov-ostrovskii}.

\subsection{Isotropic stable Lévy process}

We have seen that the boundary trace of the $(d + 1)$-dimensional Brownian motion $(X_t, Y_t)$ in the half-space is the isotropic $1$-stable Lévy process $Z_s$, and the generator of $Z_s$ coincides (up to a sign) with the Dirichlet-to-Neumann operator $K = \sqrt{-\Delta}$. What if we replace the Laplace operator by the operator $L$ introduced above? The appropriately weighted Dirichlet-to-Neumann operator becomes the fractional Laplace operator $c_\alpha (-\Delta)^{\alpha/2}$. But what happens to the boundary trace process?

Recall that
\[
 L = \Delta_x + \biggl(\frac{\partial^2}{\partial y^2} + \frac{1 - \alpha}{y} \, \frac{\partial}{\partial y} \biggr) .
\]
The operator $\tfrac12 L$ is the generator of a $(d + 1)$-dimensional diffusion process $(X_t, Y_t)$ with independent coordinates $X_t$ and $Y_t$. Furthermore, $X_t$ is the $d$-dimensional Brownian motion (with generator $\tfrac12 \Delta$), while as long as it remains positive, $Y_t$ is an Itô diffusion of the form
\[
 dY_t = dW_t + \frac{1 - \alpha}{2} \, \frac{1}{Y_t} \, dt ,
\]
where $W_t$ is some one-dimensional Brownian motion. This is a familiar stochastic differential equation: it defines the Bessel process of (fractional) dimension $2 - \alpha$. The point $0$ is a regular boundary point (so both entrance and exit from the boundary are possible) for $Y_t$ precisely when $\alpha \in (0, 2)$, and so we may require $Y_t$ to be reflected at $0$. Then there is a way to define the local time $L_t$ of $Y_t$ at $0$ (with a different scaling of approximate local times, though), and we may again define the inverse local time $T_s$ and the boundary trace process $Z_s = X(T_s)$.

Once again $Z_s$ is a Lévy process, and a straightforward scaling argument shows that $Z_s$ is self-similar with Hurst exponent $\tfrac1\alpha$. Thus, $Z_s$ is the isotropic $\alpha$-stable Lévy process, and so the corresponding generator is, up to a constant factor, the fractional Laplace operator $-(-\Delta)^{\alpha/2}$.

The result sketched above was observed and rigorously proved by Molchanov and Ostrovskii already in 1969~\cite{molchanov-ostrovskii}. Up to multiplication by a (positive or negative) constant, the generator of the boundary trace process $Z_s$ again coincides with the Dirichlet-to-Neumann operator $K$: both are equal to the fractional Laplace operator $c_\alpha (-\Delta)^{\alpha / 2}$. The density of the Lévy measure of $Z_s$ is again equal to the non-local kernel of $K$. For $\alpha \ne 1$, however, the density of the distribution of $Z_s$ is no longer equal to the Poisson kernel: we have evaluated the Fourier transform of $P_y$ and it is not of the form $\exp(-s |\xi|^\alpha)$.

The Molchanov--Ostrovskii boundary trace technique described above is a perfect probabilistic analogue of the Caffarelli--Silvestre extension technique in the theory of PDEs. It allows one to study the isotropic stable Lévy processes using the methods available for diffusions. In the remaining sections we explore to what extent we can generalise the results of Caffarelli and Silvestre on the PDE side, and of Molchanov and Ostrovskii on the probability side, and why these two approaches are in fact essentially equivalent.

\subsection*{Exercises}

\begin{problem}
\begin{enumerate}[label={\textup{(\alph*)}}]
\item
Verify that the Poisson kernel for the Caffarelli--Silvestre extension indeed satisfies the modified Laplace equation.
\item
Check that with an appropriate choice of the constant $c_{d, \alpha}$, the Poisson kernel has integral $1$.
\item
Bonus problem: Find a closed-form expression for this constant $c_{d, \alpha}$.
\end{enumerate}
\end{problem}

\begin{problem}
\begin{enumerate}[label={\textup{(\alph*)}}]
\item
Prove rigorously that the Dirichlet-to-Neumann operator $K$ for the Caffarelli--Silvestre extension is, up to a constant factor, indeed the fractional Laplace operator (given as the principal value integral).
\item
Bonus problem: Find a closed-form expression for that constant factor.
\end{enumerate}
\end{problem}

\begin{problem}
Substitute $y = z^{1 / \alpha}$ and prove that in variables $(x, z)$ the extension $u$ satisfies
\[
 \alpha^{-2} z^{2 / \alpha - 2} \Delta_x u(x, z) + \frac{\partial^2 u}{\partial z^2}(x, z) = 0 .
\]
Prove that in these coordinates the fractional Laplace operator is, up to a constant factor, again equal to the usual Dirichlet-to-Neumann operator:
\[
 K f(x) = -\frac{\partial u}{\partial z}(x, 0) ,
\]
\end{problem}

\begin{problem}
Bonus problem: Evaluate the constants in the solution of the Caffarelli--Silvestre extension problem involving Fourier transforms.
\end{problem}

\begin{problem}
\begin{enumerate}[label={\textup{(\alph*)}}]
\item
Let $Y_0 = 0$. Prove that the Bessel process $Y_t$ is self-similar with Hurst exponent $\tfrac12$.
\item
The local time of $Y_t$ at $0$ needs to be defined with a different scaling factor:
\[
 L_t = \lim_{\delta \to 0^+} \frac{L_t^{(\delta)}}{2 \delta^{2 - \alpha}} \, .
\]
Prove that $L_t$ is self-similar with Hurst exponent $\tfrac\alpha2$.
\item
Observe that the inverse local time $T_s$ is self-similar with Hurst exponent $\tfrac2\alpha$, and therefore it is the $\tfrac\alpha2$-stable subordinator.
\item
Conclude that the boundary trace process $Z_s$ is self-similar with Hurst exponent $\tfrac1\alpha$, and hence it is the isotropic $\alpha$-stable Lévy process.
\end{enumerate}
\end{problem}

%
%

\newpage
\section{Dirichlet-to-Neumann operators and Krein's strings}

\subsection{General nonsense}

At this point, our programme is hopefully clear: we study differential operators $L$ in the half-space $D = \R^d \times (0, \infty)$, and reflected diffusions $(X_t, Y_t)$ with generators $\tfrac12 L$. We would like to establish a link between the Dirichlet-to-Neumann operator $K$ corresponding to $L$, and the generator of the boundary trace of $(X_t, Y_t)$.

What class of operators can we handle? In principle, we should be able to cover essentially arbitrary elliptic second-order differential operators $L$, as long as they generate a reasonably nice diffusion process. After all, all that we need is the Poisson kernel on the PDE side, and the local time on the probability side. A general treatment requires precise assumptions on the diffusion, its reflection mechanism, and the boundary operator.

Nevertheless, we will not develop that theory here. (If you are interested in doing so, however, do let me know!) The reasons for that are two-fold. First, we have in mind applications to non-local PDEs in $\R^d$, and so it is important for us to understand the class of operators $K$ that we can handle. Second, on the probability side we do want to work with the usual definition of the local time.

\subsection{Elliptic operators}

In what follows we consider elliptic operators of the form
\[
 L = a(y) \Delta_x + b(y) \, \frac{\partial}{\partial y} + c(y) \, \frac{\partial^2}{\partial y^2} \, ,
\]
where $a$, $b$ and $c$ are appropriately regular coefficients. In other words, we require $L$ to be invariant under isometric transformations (rotations, translations, reflections) of variable~$x$.

Since we are only concerned with functions harmonic with respect to $L$, we may freely multiply $L$ by a function of $y$. In probabilistic terms, this corresponds to a random time change, which preserves the boundary trace after the corresponding adjustment of its local-time parametrisation.

Similarly, we can transform the coordinate $y$ in a strictly increasing (possibly nonlinear) way. This change of scale preserves both the Dirichlet-to-Neumann operator and the boundary trace process, as long as we take the change of scale into account and appropriately adjust the denominators in the definitions of the Dirichlet-to-Neumann operator $K$ and the local time~$L_t$.

By combining the above two operations, and ignoring completely regularity issues, we can cast an arbitrary operator $L$ described above into a form that is convenient for our needs, essentially reducing the number of coefficients depending on $y$ from three ($a$, $b$ and~$c$) to one. A natural choice from the point of view of applications in probability is the operator
\[
 L = \Delta_x + b(y) \, \frac{\partial}{\partial y} + \frac{\partial^2}{\partial y^2} \, ,
\]
which is (up to a constant factor $\tfrac12$) the generator of a particularly simple Itô diffusion, the Brownian motion with a state-dependent drift. In PDEs, one would likely prefer the same operator written in divergence form
\[
 L = \nabla_{x, y} \cdot (\beta(y) \nabla_{x, y}) .
\]
The former operator can be obtained from the latter by dividing by $\beta(y)$ and setting $b(y) = \beta'(y) / \beta(y)$; going in the opposite direction is equally simple. These two forms are the ones that we used in the Caffarelli--Silvestre extension technique and in the Molchanov--Ostrovskii boundary trace technique.

When it comes to actually proving things, it is, however, much more convenient to work with the operators
\[
 L = a(y) \Delta_x + \frac{\partial^2}{\partial y^2} \, ,
\]
where $a$ is a non-negative locally integrable function. Additionally, with this form of $L$ the Dirichlet-to-Neumann operator is defined in the usual way, with just $y$ in the denominator. Finally, the above form of $L$ allows for a generalisation which makes our results complete.

The last sentence is to be understood as follows: if we agree to include operators on $\R^d \times [0, R)$ with a possibly finite $R \in (0, \infty]$ (and with an additional zero boundary condition on $\R^d \times \{R\}$ whenever necessary), and if we allow the coefficient $a$ to be a locally finite non-negative measure on $[0, R)$, then we obtain a perhaps surprising `if and only if' result. To keep things simple, though, with few exceptions we will always pretend that $a$ is a function and $R = \infty$. That said, you may verify that the general case requires only minor modifications indicated below.

\subsection{Krein's spectral theory of strings}

After this very long introduction, the fun part begins: we come to the heart of the PDE section of these notes.

Applying the Fourier transform to both sides of the equation $L u = 0$, with $L$ defined above, leads to the ODE
\[
 \frac{\partial^2 \fourier u}{\partial y^2}(\xi, y) = |\xi|^2 a(y) \fourier u(\xi, y) ,
\]
with initial condition $\fourier u(\xi, 0) = \fourier f(\xi)$, and with $\fourier u(\xi, y)$ bounded with respect to $y$.

If we substitute $\lambda = |\xi|^2$ and $\ph(y) = \fourier u(\xi, y)$, we arrive at the eigenvalue problem for a second-order ordinary differential operator, studied by Krein:
\[
 \ph''(y) = \lambda a(y) \ph(y) .
\]
If $\ph_\lambda$ is a bounded solution satisfying $\ph_\lambda(0) = 1$, then the extension $u$ is given by
\[
 \fourier u(\xi, y) = \ph_\lambda(y) \fourier f(\xi) ,
\]
where $\lambda = |\xi|^2$.

It is a routine application of Picard's iteration to prove that for every complex $\lambda$, there exists a solution $\ph$ of Krein's ODE with any prescribed $\ph(0)$ and $\ph'(0)$. Let us denote by $\ph_\dir$ the `Dirichlet' solution with $\ph_\dir(0) = 0$, $\ph_\dir'(0) = 1$, and let $\ph_\neu$ be the `Neumann' solution with $\ph_\neu(0) = 1$, $\ph_\neu'(0) = 0$. Our goal is to construct a solution $\ph = \ph_\neu - \mu \ph_\dir$ satisfying the initial condition $\ph(0) = 1$ and which is bounded on $[0, \infty)$. The last condition is a kind of a `terminal' condition at $\infty$, so essentially we are solving a variant of the Sturm--Liouville eigenvalue problem.

When $\lambda \ge 0$, both fundamental solutions $\ph_\dir$ and $\ph_\neu$ are non-decreasing, and we take
\[
 \mu = \lim_{y \to \infty} \frac{\ph_\neu(y)}{\ph_\dir(y)} \, .
\]
It is fairly easy to see that $\ph_\neu(y) / \ph_\dir(y)$ is decreasing:
\[
 \biggl(\frac{\ph_\neu}{\ph_\dir}\biggr)^{\!\prime} = \frac{\ph_\dir \ph_\neu' - \ph_\neu \ph_\dir'}{(\ph_\dir)^2}
\]
and since the Wronskian $\ph_\dir \ph_\neu' - \ph_\neu \ph_\dir'$ has constant sign, the right-hand side is negative. This property implies that $\mu$ is well-defined, and additionally $\ph_\neu / \ph_\dir \ge \mu$, that is, $\ph_\neu - \mu \ph_\dir \ge 0$. A slightly more careful analysis reveals that $\ph_\neu - \mu \ph_\dir$ is in fact non-increasing. The argument is a bit awkward, we pass to the limit as $y^* \to \infty$ after making the following point: the function $\ph_\neu - (\ph_\neu(y^*) / \ph_\dir(y^*)) \ph_\dir$ is positive up to $y^*$ and then negative (by monotonicity of $\ph_\neu / \ph_\dir$), and so also convex up to $y^*$ and then concave (because it is a solution of Krein's ODE), and this is only possible if it is decreasing.

Since non-increasing and non-negative functions on $[0, \infty)$ are bounded, we see that
\[
 \ph = \ph_\neu - \mu \ph_\dir
\]
is the solution we have been looking for. Clearly, we have $\ph'(0) = -\mu$. What we have proved so far is summarised in the following statement.

\begin{proposition}
Given $\lambda \ge 0$, there is a unique bounded solution $\ph_\lambda$ of Krein's ODE
\[
 \ph''(y) = \lambda a(y) \ph(y)
\]
satisfying $\ph_\lambda(0) = 1$.
\end{proposition}

Let us comment briefly on the general case. When $a$ is a measure, Krein's ODE should be understood in the integral form
\[
 \ph(y) = \ph(0) + \ph'(0) y + \int_0^y \biggl(\int_{[0, y')} \lambda \ph(y'') a(dy'') \biggr) dy' .
\]
In particular, if $a(\{0\}) > 0$, then $\ph'(0)$ has a rather unexpected, and very unfortunate, meaning: the right-hand side derivative of $\ph$ at $0$ is equal to $\ph'(0) + \lambda a(\{0\}) \ph(0)$ (for this reason some authors tend to use $\ph'(0^-)$ instead). If $R$ is finite and $\int_{[0, R)} (R - y) a(dy)$ is finite, then every solution $\ph$ is bounded, and the one we are interested in satisfies the Dirichlet boundary condition $\ph(R) = 0$. Otherwise, the argument is exactly the same. These remarks implicitly apply to the proposition given above, as well as to all results that follow.

If $R$ is finite and $\int_{[0, R)} (R - y) a(dy)$ is infinite and $\lambda > 0$, then bounded solutions automatically satisfy the Dirichlet condition $\ph(R) = 0$, which we also impose for $\lambda = 0$. In particular, we have $\ph_0(y) = 1$ and $\mu(0) = 0$ when $R = \infty$, while $\ph_0(y) = 1 - \frac{y}{R}$ and $\mu(0) = \frac{1}{R} > 0$ otherwise.

The measure $a$ can be understood as the distribution of mass of a (perfectly elastic) vibrating string. For this reason, Krein used the word string for the distribution function of $a$, that is, $y \mapsto a([0, y))$, and the study of the spectral problem discussed above is widely known as Krein's theory of strings.

We will shortly see that the Dirichlet-to-Neumann operator $K$ corresponding to the elliptic operator $L$ is a multiplication operator in the Fourier space, with symbol $\mu(|\xi|^2)$, where $\mu(\lambda) = -\ph_\lambda'(0)$ (this is the same number $\mu$ as the one constructed above). So what do we know about the function $\mu$? The answer is provided by the following beautiful and deep result proved by Krein in 1952~\cite{krein}.

\begin{theorem}
\begin{enumerate}[label={\textup{(\alph*)}}]
\item
If $\mu(\lambda) = -\ph_\lambda'(0)$, then $\mu$ is a complete Bernstein function. By this we mean that $\mu : [0, \infty) \to [0, \infty)$ is continuous and extends to a holomorphic function on $\C \setminus (-\infty, 0]$ that maps the upper and lower complex half-planes into their respective closures. Equivalently:
\[
 \mu(\lambda) = p + q \lambda + \int_{(0, \infty)} \frac{\lambda}{\lambda + s} \, m(ds)
\]
for some $p, q \ge 0$ and a non-negative measure $m$ satisfying $\int_{(0, \infty)} (1 + s)^{-1} m(ds) < \infty$.
\item
Every complete Bernstein function can be represented in the way described above, for a unique coefficient $a$.
\end{enumerate}
\end{theorem}

We stress that when we say `for a unique coefficient $a$', what we really mean is: for a unique $R \in (0, \infty]$ and a unique non-negative locally finite measure $a(dy)$ on $[0, R)$.

Krein's theorem is not constructive: there is no simple explicit reconstruction formula. Given the coefficient $a$, if one is lucky, one can solve the corresponding ODE and get an expression for the complete Bernstein function $\mu$. However, finding a closed-form inverse is generally difficult unless $\mu$ has a special form. (For example, if $\mu$ is a rational function, then the corresponding coefficient $a$ can be reconstructed using an appropriate continued fraction expansion. In this case, however, $a$ turns out to be a purely atomic measure.)

Finally, we remark that $1 / \mu$ is often called the spectral, Weyl or characteristic function.

\subsection{Dirichlet-to-Neumann operators}

Our construction of the solution of Krein's spectral problem given in the first proposition above immediately leads to the following statement.

\begin{corollary}
Let us write $\ph(\lambda, y)$ instead of $\ph_\lambda(y)$ for the bounded solution of Krein's ODE satisfying $\ph_\lambda(0) = 1$. Given a square integrable function $f$ on $\R^d$, define $u$ in terms of the Fourier transform:
\[
 \fourier u(\xi, y) = \ph(|\xi|^2, y) \fourier f(\xi) .
\]
Then, in an appropriate $L^2$ sense, $u$ is a solution of the PDE
\[
 a(y) \Delta_x u(x, y) + \frac{\partial^2 u}{\partial y^2}(x, y) = 0 ,
\]
and the $L^2$ norm of $u(\cdot, y)$ does not exceed the $L^2$ norm of $f$ for every $y > 0$.
\end{corollary}

Of course, the `appropriate $L^2$ sense' above means in particular that the Fourier transform of $u$ is the `usual' solution of Krein's ODE for almost every parameter $\xi$, and at the same time $y \mapsto u(\cdot, y)$ is a bounded and continuous $L^2$-valued function. More can be said about~$u$: for example, it satisfies the PDE in the sense of distributions; but we will not pursue that direction.

The above corollary in particular means that the Poisson kernel $P_y$ for the operator~$L$ in the upper half-space has Fourier transform
\[
 \fourier P_y(\xi) = \ph(|\xi|^2, y) .
\]
In the next section, using probabilistic tools, we will easily verify that $P_y$ is a non-negative measure of total mass at most one. Interestingly, a purely analytical proof of this fact seems to require much more effort.

Since $\ph_\lambda$ is a convex non-increasing function, the difference quotient of $-\ph_\lambda$ on the interval $[0, y]$ is positive and it increases to $\mu(\lambda) = -\ph_\lambda'(0)$ as $y$ decreases to $0$. Thus, Lebesgue's dominated convergence theorem implies that the Dirichlet-to-Neumann operator $K$ corresponding to the elliptic operator $L$ is given by

For a function $a$ (or a measure without an atom at $0$), $\ph_\lambda$ is convex, and the quantity $(1-\ph_\lambda(y))/y$ is non-negative and it increases to $\mu(\lambda) = -\ph_\lambda'(0)$ as $y$ decreases to $0$. Thus, Lebesgue's dominated convergence theorem implies that the Dirichlet-to-Neumann operator $K$ corresponding to the elliptic operator $L$ is given by
\[
 \fourier K f(\xi) = -\frac{\partial \fourier u}{\partial y}(\xi, 0) = -\frac{\partial \ph}{\partial y}(|\xi|^2, 0) \fourier f(\xi) = \mu(|\xi|^2) \fourier f(\xi)
\]
as long as the function $\mu(|\xi|^2) \fourier f(\xi)$ is square integrable. We have thus essentially proved the first two parts of the following result, given in my joint paper with Jacek Mucha from 2018~\cite{kwasnicki-mucha}.

\begin{proposition}
\begin{enumerate}[label={\textup{(\alph*)}}]
\item
Dirichlet-to-Neumann operators $K$ corresponding to elliptic operators of the form
\[
 L = a(y) \Delta_x + b(y) \, \frac{\partial}{\partial y} + c(y) \, \frac{\partial^2}{\partial y^2}
\]
are complete Bernstein functions of $-\Delta$. That is, there is a complete Bernstein function $\mu$ such that
\[
 \fourier K f(\xi) = \mu(|\xi|^2) \fourier f(\xi)
\]
whenever $f$ and the right-hand side are square integrable.
\item
Every complete Bernstein function of $-\Delta$ can be represented in the way described above, for a unique operator of the form
\[
 L = a(y) \Delta_x + \frac{\partial^2}{\partial y^2} \, .
\]
\item
For $L$ as in the previous item, the following Dirichlet's principle holds:
\[
 \begin{aligned}
  & \int_{\R^d} f(x) K f(x) dx \\
  & \hspace*{3em} = \min\biggl\{\int_{\R^d} \int_0^\infty \biggl(a(y) |\nabla_x v(x, y)|^2 + \biggl(\frac{\partial v}{\partial y}(x, y)\biggr)^{\!\!2\,}\biggr) dy dx\biggr\} ,
 \end{aligned}
\]
where the minimum is taken over all functions $v$ which are suitable differentiable, are suitably bounded as $y \to \infty$, and have boundary values $f$.
\end{enumerate}
\end{proposition}

Complete Bernstein functions have been studied a lot, see Chapters~6 and~7 in the excellent book~\cite{schilling-song-vondracek} by Schilling, Song and Vondraček. Complete Bernstein functions of $-\Delta$ play an important role in the theory of non-local PDEs and in the theory of Lévy processes.

\subsection*{Exercises}

\begin{problem}
Prove that the Dirichlet-to-Neumann operator corresponding to the operator
\[
 L = \ind_{(0, 1)}(y) \Delta_x + \frac{\partial^2}{\partial y^2}
\]
in $\R^d \times (0, \infty)$, is $K = \sqrt{-\Delta} \tanh \sqrt{-\Delta}$, in the sense that the complete Bernstein function arising in the corresponding Krein's spectral problem is $\mu(\lambda) = \sqrt{\lambda} \tanh \sqrt{\lambda}$. This operator $K$ arises naturally in the theory of linear water waves.
\end{problem}

\begin{problem}
Prove that the Dirichlet-to-Neumann operator corresponding to the operator
\[
 L = \Delta_{x, y} \qquad \text{in $\R^d \times (0, 1)$,}
\]
with an implicit Dirichlet condition $u(x, 1) = 0$, is $K = \sqrt{-\Delta} \coth \sqrt{-\Delta}$.
\end{problem}

\begin{problem}
Prove that the Dirichlet-to-Neumann operator corresponding to the operator
\[
 L = \frac{\partial^2}{\partial y^2}
\]
in $\R^d \times (0, \infty)$ is the zero operator $K = 0$.
\end{problem}

\begin{problem}
Prove that the Dirichlet-to-Neumann operator corresponding to the operator
\[
 L = \frac{\partial^2}{\partial y^2} \qquad \text{in $\R^d \times (0, 1)$,}
\]
with an implicit Dirichlet condition $u(x, 1) = 0$, is the identity operator $K = 1$.
\end{problem}

\begin{problem}
Prove that the Dirichlet-to-Neumann operator corresponding to the operator
\[
 L = \delta_1(dy) \Delta_x + \frac{\partial^2}{\partial y^2}
\]
in $\R^d \times (0, \infty)$ is $K = -\Delta (1 - \Delta)^{-1}$.
\end{problem}

\begin{problem}
Prove that the Dirichlet-to-Neumann operator corresponding to the operator
\[
 L = \delta_0(dy) \Delta_x + \frac{\partial^2}{\partial y^2}
\]
in $\R^d \times (0, \infty)$ is, up to a sign, the Laplace operator $K = -\Delta$.
\end{problem}

\begin{problem}
Prove that the Dirichlet-to-Neumann operator corresponding to the operator
\[
 L = \frac{1}{(1 + 2 y)^2} \, \Delta_x + \frac{\partial^2}{\partial y^2}
\]
in $\R^d \times (0, \infty)$ is the quasi-relativistic operator $K = \sqrt{-\Delta + 1} - 1$.
\end{problem}

\begin{problem}
Prove that the Dirichlet-to-Neumann operator corresponding to the operator
\[
 L = \frac{1}{(1 - 2 y)^2} \, \Delta_x + \frac{\partial^2}{\partial y^2} \qquad \text{in $\R^d \times (0, \tfrac12)$}
\]
is $K = \sqrt{-\Delta + 1} + 1$. (Note that the Dirichlet condition $u(x, \tfrac12) = 0$ is automatically satisfied by bounded solutions.)
\end{problem}

\begin{problem}
Bonus problem: Find an elliptic operator $L$ such that the corresponding Dirichlet-to-Neumann operator is $K = \sqrt{-\Delta + 1}$.
\end{problem}

%
%

\newpage
\section{Boundary traces and Dirichlet-to-Neumann operators}

\subsection{Intuitive argument}

We now study essentially the same problem from the probabilistic perspective: our goal is to find a link between boundary traces of diffusions and Dirichlet-to-Neumann operators $K$ corresponding to their generators $L$. This may appear nearly obvious at an intuitive level, but a direct approach does not seem feasible unless some regularity of the coefficients of $L$ near the boundary is assumed. For this reason, we tackle the problem in a more analytical fashion. Before we proceed, however, let us make one remark, and then try a different, seemingly more probabilistic approach.

We consider a diffusion process $(X_t, Y_t)$ in $\R^d \times [0, \infty)$, reflected at the boundary. The boundary trace process is given by $Z_s = X(T_s)$, where $T_s$ is the inverse local time of $Y_t$ at $0$. It is clear that the jumps of $Z_s$ correspond to excursions of $(X_t, Y_t)$ away from the boundary. One can prove that if $(X_0, Y_0) = (x, y)$ and if $\tau_0$ is the hitting time of $0$ for $Y_t$, then the law of $X(\tau_0)-x$ is the Poisson kernel $P_y$. Thus, as $y \to 0^+$, appropriately normalised Poisson kernels are expected to converge to the distribution of the endpoint for excursions started at $(x, 0)$, which we already identified with the Lévy measure (that is, the intensity of jumps) of $Z_s$. On the other hand, appropriately normalised Poisson kernels $P_y(x - x')$ converge to the non-local kernel of the Dirichlet-to-Neumann operator. This may convince us that indeed these two kernels are the same thing. A rigorous argument of that kind, however, seems to be out of reach.

\subsection{Generalised diffusions}

In the previous section we studied elliptic operators of the form
\[
 L = a(y) \Delta_x + \frac{\partial^2}{\partial y^2} \, ,
\]
and now we would like to study diffusion processes generated by these operators (or, strictly speaking, by $\tfrac12 L$). Recall, however, that we have much flexibility here: we can change the above operator $L$ to a different, equivalent form using change of scale and multiplication by a given function. Therefore, it seems reasonable to consider instead an equivalent operator
\[
 L = \Delta_x + \frac{1}{a(y)} \frac{\partial^2}{\partial y^2} \, ,
\]
or perhaps the operator which corresponds to the Brownian motion with a state-de\-pen\-dent drift:
\[
 L = \Delta_x + b(y) \, \frac{\partial}{\partial y} + \frac{\partial^2}{\partial y^2} \, .
\]
These operators have two separate parts, one of which is simply $\Delta_x$, and the other one is given purely in terms of variable $y$. Thus, the diffusion $(X_t, Y_t)$ with generator $\tfrac12 L$ has independent components, and the first one, $X_t$, is simply the Brownian motion in $\R^d$. This independence immediately tells us what we get: the trace process $Z_s = X(T_s)$ is the subordinate Brownian motion, determined by the subordinator $T_s$ obtained as the inverse local time of $Y_t$ at $0$. Therefore, it remains to determine what kind of subordinators $T_s$ can be obtained in this way.

This question was raised by Itô and McKean in 1965 in their book~\cite{ito-mckean}. While a complete characterisation for Itô diffusions is not known, an answer for a broader class of generalised diffusions, or gap diffusions, was given independently by Knight in 1981~\cite{knight} and by Kotani and Watanabe in 1982~\cite{kotani-watanabe}: the subordinator $T_s$ can be obtained as the inverse local time of a reflected generalised diffusion (allowing killing and degenerate cases) if and only if its characteristic (Laplace) exponent is a complete Bernstein function.

As you have certainly guessed, the generators of generalised diffusions are precisely operators of the form
\[
 \frac{1}{a(y)} \frac{\partial^2}{\partial y^2} \, ,
\]
on an infinite or finite interval $[0, R)$, where in fact we allow $a$ to be a locally finite measure. And it will be no surprise to you that the main tool employed in this result is Krein's spectral theory of strings.

While the above argument is in some sense satisfactory, arguably it provides little insight into the relation between the Dirichlet-to-Neumann operator and the generator of the boundary trace process. These two turn out to coincide, but except for the fact that the same tool is used as a black box in the proofs, the two arguments seem very loosely related. It is therefore difficult to apply the same approach if we allow for an extra term in $L$ of the form $\sum_{j = 1}^d b_j(y) \frac{\partial^2}{\partial x_j \partial y}$.

Let us make one more remark here. The name gap diffusion refers to the fact that if $a(y) = 0$ on some maximal interval, then the process $Y_t$ jumps between the endpoints of that interval and it never enters its interior. If this phenomenon occurs in the vicinity of the boundary point $0$, then our intuition with excursions and limits of the Poisson kernel described in the beginning of this section breaks down.

\subsection{MAP}

Quite surprisingly, one gets far better results by studying the same operator as in the previous section:
\[
 L = a(y) \Delta_x + \frac{\partial^2}{\partial y^2} \, .
\]
The operator $\tfrac12 L$ is the generator of a diffusion $(X_t, Y_t)$ for which $X_t$ is in general not the Brownian motion, and $X_t$ and $Y_t$ need not be independent. On the other hand, however, $Y_t$ is simply the reflected Brownian motion in $[0, \infty)$, which regulates the speed of the evolution of the, otherwise Brownian, evolution of $X_t$.

This means that the pair $(Y_t, X_t)$ is a Markov additive process, abbreviated as MAP, with $Y_t$ being the regulator term, and $X_t$ the additive term. Fortunately, in this section we will not need any tools from the theory of MAPs. What follows is essentially identical to the calculation carried out in a slightly different context in my article from 2023~\cite{kwasnicki-2}.

What we need is a process $(X_t, Y_t)$ where $Y_t$ is the reflected Brownian motion in $[0, \infty)$, and, conditionally on the path of $Y_t$, $X_t$ is the time-changed Brownian motion characterised by $d\langle X\rangle_t = a(Y_t) \id dt$. We can construct $X_t$ using an independent Brownian motion $W_t$ in $\R^d$ by setting 
\[
 dX_t = \sqrt{a(Y_t)} \, dW_t .
\]
Clearly, the quadratic variation of $X_t$ satisfies $d\langle X\rangle_t = a(Y_t) \id dt$, as desired. It turns out, however, that there is a better way to go.

Let us again start with the same pair $(W_t, Y_t)$, and denote by
\[
 A_t = \int_0^t a(Y_s) ds
\]
the non-decreasing additive functional of $Y_t$ (with Revuz measure $a(y) dy$ if one is familiar with this concept). Now we set
\[
 X_t = W(A_t) .
\]
Then $X_t$ is a continuous local martingale with quadratic variation $\langle X\rangle_t = A_t \id$, and so the pair $(X_t, Y_t)$ is exactly what we need.

It is the right moment to mention that when $a$ is a measure rather than a function, then the above definitions gets slightly more complicated: in order to properly define $A_t$, we integrate with respect to $a(dy)$ the family of local times of $Y_t$ at all points $y \ge 0$, appropriately normalised. But this is a standard procedure in this business (in sharp contrast to the problems that occur with our first approach, $dX_t = \sqrt{a(Y_t)} \, dW_t$).

\subsection{Itô calculus}
\label{sec:ito}

The main advantage of the above construction is that it allows us to calculate various quantities related to $(X_t, Y_t)$ in an effective way, and link them with the objects that we used in the previous section. In order to do so, we get rid of $X_t$ and the underlying Brownian motion $W_t$ by considering a probabilistic counterpart of the Fourier transform with respect to variable $x$. Namely, for a given value of the Fourier parameter $\xi$ we evaluate the `partial' expectation of
\[
 \exp(-i \xi X_t)
\]
with respect to the process $W_t$, while leaving the process $Y_t$ random. This is formally written as the conditional expectation with respect to the $\sigma$-algebra generated by the entire path of $Y_t$, but we prefer to use a simpler (and more common sense) notation $\ex_W$. Also, here and below, with no loss of generality, we assume that $(X_0, Y_0) = (0, 0)$. Thus, we define
\[
 \Xi_t = \ex_W [\exp(-i \xi X_t)] .
\]
Since $X_t = W(A_t)$ and $A_t$ is a functional of $Y_t$, which is independent of $W_t$, we simply have
\[
 \Xi_t = \ex_W [\exp(-i \xi W(A_t))] = \exp(-\tfrac12 |\xi|^2 A_t) .
\]
In other words, $\Xi_t$ is a non-increasing process started at $\Xi_0 = 1$, satisfying
\[
 \begin{aligned}
  d\Xi_t & = -\tfrac12 |\xi|^2 \Xi_t dA_t \\
  & = -\tfrac12 |\xi|^2 \Xi_t a(Y_t) dt .
 \end{aligned}
\]

Now we only have one source of randomness left: the reflected Brownian motion $Y_t$ in $[0, \infty)$. The term $|\xi|^2 a(y)$ on the right-hand side of the expression for $d\Xi_t$ has already appeared above, in the right-hand side of Krein's ODE
\[
 \ph''(y) = |\xi|^2 a(y) \ph(y)
\]
studied in the previous section. So let us denote by $\ph$ the bounded solution which satisfies $\ph(0) = 1$, denoted by $\ph(|\xi|^2, \cdot)$ in the previous section, and let us inspect the process
\[
 \Phi_t = \ph(Y_t) .
\]
By Itô's lemma,
\[
 \begin{aligned}
  d\Phi_t & = \ph'(Y_t) dY_t + \tfrac12 \ph''(Y_t) dt \\
  & = \ph'(Y_t) dY_t + \tfrac12 |\xi|^2 a(Y_t) \Phi_t dt .
 \end{aligned}
\]
Therefore, the bounded variation terms in $\Xi_t$ and $\Phi_t$ will cancel out if we multiply these two processes. Indeed: using the integration by parts formula for Itô integrals:
\[
 d(\Xi \Phi)_t = \Xi_t d\Phi_t + \Phi_t d\Xi_t + d[\Xi, \Phi]_t ,
\]
and the fact that $\Xi$ is of bounded variation, we arrive at
\[
 \begin{aligned}
  d(\Xi \Phi)_t & = \Xi_t \cdot \bigl(\ph'(Y_t) dY_t + \tfrac12 |\xi|^2 a(Y_t) \Phi_t dt\bigr) + \Phi_t \cdot \bigl(-\tfrac12 |\xi|^2 \Xi_t a(Y_t) dt\bigr) + 0 \\
  & = \Xi_t \ph'(Y_t) dY_t .
 \end{aligned}
\]
So it turns out that $\Xi_t \Phi_t$ is almost a martingale: $Y_t$ is the reflected Brownian motion in $[0, \infty)$, so $dY_t$ is essentially the differential of the Brownian motion, except when $Y_t = 0$. More precisely, if $L_t$ is the local time of $Y_t$ at zero, then, by Tanaka's formula, $\tilde Y_t = Y_t - L_t$ is the (unrestricted) Brownian motion. Therefore, the Doob--Meyer decomposition of $\Xi_t \Phi_t$ is given by
\[
 d(\Xi \Phi)_t = \Xi_t \ph'(Y_t) d\tilde Y_t + \Xi_t \ph'(Y_t) dL_t ,
\]
the former term being the differential of a local martingale, and the latter one the differential of a bounded variation process. The local martingale stopped at $T_s$ turns out to be uniformly integrable, and so, integrating from $0$ to $T_s$ and applying the optional stopping theorem, its expectation vanishes. We arrive at
\[
 \ex_Y [\Xi(T_s) \Phi(T_s)] = \ex_Y [\Xi_0 \Phi_0] + \ex_Y \biggl[ \int_0^{T_s} \Xi_t \ph'(Y_t) dL_t \biggr]
\]
(we write $\ex_Y$ to emphasise that we have already integrated over the law of $W_t$).

To take a breath before we move on, let us appreciate the fact that we have just used four significant results from the stochastic analysis toolkit: Itô's lemma, integration by parts, Tanaka's formula and optional stopping theorem, all in a single paragraph.

We are almost done. Since we assumed that $X_0 = 0$ and $Y_0 = 0$, we have $\Xi_0 = \ex_W \exp(-i \xi X_0) = 1$, and $\Phi_0 = \ph(Y_0) = 1$. Substituting $L_t = r$ in the integral on the right-hand side (which is a Lebesgue--Stieltjes integral, not an Itô one), and recalling that $\Phi_t = \ph(Y_t)$, we find that
\[
 \ex_Y [\Xi(T_s) \ph(Y(T_s))] = 1 + \ex_Y \biggl[ \int_0^s \Xi(T_r) \ph'(Y(T_r)) dr \biggr] .
\]
However, $\Xi_t = \ex_W [\exp(-i \xi X_t)]$, $T_s$ and $W$ are independent, and so $\Xi(T_s) = \ex_W [\exp(-i \xi X(T_s))]$. Furthermore, $X(T_s) = Z_s$ is the boundary trace process, and $Y(T_s) = 0$ by the definition of the inverse local time $T_s$. Finally, we have $\ph(0) = 1$ and $\ph'(0) = -\mu(|\xi|^2)$, where $\mu$ is a complete Bernstein function from Krein's spectral problem. Therefore, our equation reduces to
\[
 \ex_Y [\ex_W [\exp(-i \xi Z_s)] \cdot 1] = 1 + \ex_Y \biggl[ \int_0^s \ex_W [\exp(-i \xi Z_r)] \cdot (-\mu(|\xi|^2)) dr \biggr] .
\]
Using Fubini's theorem and the fact that $\ex_Y \ex_W = \ex$, we conclude that
\[
 \ex [\exp(-i \xi Z_s)] = 1 - \mu(|\xi|^2) \int_0^s \ex [\exp(-i \xi Z_r)] dr .
\]
The above equation is the integral form of a simple first-order ODE for $\psi(s) = \ex [\exp(-i \xi Z_s)]$, namely,
\[
 \psi'(s) = -\mu(|\xi|^2) \psi(s)
\]
with initial condition $\psi(0) = 1$. A solution is an appropriate exponential function, and we conclude that
\[
 \ex [\exp(-i \xi Z_s)] = \exp(-s \mu(|\xi|^2)) .
\]
But hey, this is the Lévy--Khintchine formula for the Lévy process $Z_s$! We have just proved that $\mu(|\xi|^2)$ is the characteristic exponent of the boundary trace process $Z_s$, and this is exactly what we wanted.

I recommend you to follow the above argument once more, bearing in mind the final result. Our conclusion may have appeared in an unexpected way from seemingly unrelated calculations. On the second reading, however, the derivation of the Lévy--Khintchine formula for $Z_s$ should appear much clearer and way more natural.

\subsection*{Exercises}

\begin{problem}
Follow the recommendation from the last paragraph: read again Section~\ref{sec:ito} and try to figure out why it is in a sense natural to consider the process $\Phi_t = \ph(Y_t)$ and multiply it by~$\Xi_t$.
\end{problem}

\begin{problem}
\begin{enumerate}[label={\textup{(\alph*)}}]
\item
Redo the calculation from Section~\ref{sec:ito} with an arbitrary starting point $(X_0, Y_0) = (x, y)$, and stop the process $\Xi_t \Phi_t$ at the first hitting time of $0$ for $Y_t$ (which just happens to agree with $T_0$) rather than at $T_s$.
\item
Prove that the characteristic function of $X(T_0)$ is equal to $e^{-i \xi x} \ph(|\xi|^2, y)$, and so the law of $X(T_0) - x$ is the Poisson kernel $P_y$ discussed in the previous section.
\item
Show that if $f$ is non-negative, then also its extension $u$, introduced in the previous section, is non-negative.
\item
Prove that if $f$ is a compactly supported continuous function, then its extension $u$ is continuous.
\end{enumerate}
\end{problem}

\begin{problem}
Let $(X_n, Y_n)$ be a $(d + 1)$-dimensional simple random walk on the lattice $\Z^d \times \{0, 1, \ldots\}$, reflected on the boundary. (There are many ways to actually implement the reflection mechanism and you are welcome to choose your favourite one.) Define the local time $L_n$ of $Y_n$ at $0$ in the obvious way:
\[
 L_n = \sum_{j = 0}^n \ind_{\{0\}}(Y_j) ,
\]
and the inverse local time $T_k$ in an equally natural manner:
\[
 T_k = \min\{n\ge0:L_n\ge k+1\},\qquad k=0,1,\ldots,
\]
Assume $Y_0=0$, so that $T_0=0$. Find as much information as possible about the boundary trace process $Z_k = X(T_k)$.
\end{problem}

%
%

\newpage
\section{Extensions}

\subsection{More general symmetric operators and processes}

It is fairly clear that the Laplace operator $-\Delta_x$ in the Caffarelli--Silvestre extension technique can be replaced by an arbitrary non-negative self-adjoint operator $L_x$ on a Hilbert space $H$. One just replaces the Fourier transform $\fourier$ by the spectral resolution of $L_x$ to turn the extension problem into the same Krein's ODE. While the idea is clear, the details may be overwhelming if $L_x$ is, say, a differential operator, unless one agrees to work in the context of $H$-valued functions $y \mapsto u(\cdot, y)$. We refer to the work by Hauer and Lee~\cite{hauer-lee} for further discussion.

On the probability side, the above extension corresponds to replacing the Brownian motion $X_t$ by a more general (diffusion) process. This was carried out, with a slightly different approach, in a paper by Assing and Herman from 2021~\cite{assing-herman}.

\subsection{Non-symmetric operators and diffusions}

As suggested in the previous section, the evaluation of the Lévy--Khintchine exponent of the boundary trace process actually works for a broader class of operators $L$. We still need to assume that the coefficients only depend on $y$, so that $L$ is invariant under translations, but invariance under rotations is not really important. There is no reason to expect any reasonable class of boundary trace processes that can be obtained in this way in higher dimensions, but for $d = 1$ one gets the following, perhaps unexpected, result. We state it simultaneously in two flavours: a PDE variant and its probabilistic reformulation.

\begin{theorem}
\begin{enumerate}[label={\textup{(\alph*)}}]
\item Dirichlet-to-Neumann operators $K$ corresponding to elliptic second-order differential operators on $\R \times (0, \infty)$ of the form
\[
 L = a(y) \frac{\partial^2}{\partial x^2} + b(y) \frac{\partial^2}{\partial x \partial y} + c(y) \frac{\partial^2}{\partial y^2} + d(y) \frac{\partial}{\partial x} + e(y) \frac{\partial}{\partial y}
\]
have completely monotone kernels: they are given by
\[
 \begin{aligned}
  K f(x) & = -\alpha f''(x) + \beta f'(x) + \gamma f(x) \\
  & \hspace*{3em} + \int_{-\infty}^\infty \bigl( f(x) - f(x + z) + z f'(x) \ind_{(-1, 1)}(z) \bigr) J(z) dz ,
 \end{aligned}
\]
where $\alpha, \gamma \ge 0$, $\beta \in \R$, $\min\{1, x^2\} J(x)$ is integrable, and both $J(x)$ and $J(-x)$ are completely monotone functions of $x > 0$.
\item The correspondence between $L$ and $K$ is onto, and if we restrict our attention to the case $d(y) = e(y) = 0$, then it is additionally one-to-one.
\item The trace of a diffusion process on $\R \times [0, \infty)$, which is reflected at the boundary and which is invariant under horizontal translations, can be identified with a Lévy process with completely monotone jumps (that is, with Lévy measure $J(x) dx$ for a function $J$ as in the first item).
\item Conversely, every Lévy process with completely monotone jumps is the boundary trace process for a uniquely determined generalised diffusion in $\R \times [0, \infty)$, invariant under horizontal translations, which is additionally a martingale and the vertical component of which is the reflected Brownian motion.
\end{enumerate}
\end{theorem}

The first two statements are given in my paper from 2022~\cite{kwasnicki-1}, the other two in a paper of mine from 2023~\cite{kwasnicki-2}. The proof requires some more advanced tools. Most notably, it hinges on an extension of Krein's theory developed by Eckhardt and Kostenko in 2016~\cite{eckhardt-kostenko} (and built upon the results proved by Louis de Branges in 1960--62~\cite{debranges-1,debranges-2,debranges-3,debranges-4,debranges}).

It is difficult to predict whether the above theorems will find applications, but working out their proofs brought me a lot of joy and satisfaction. And it finally made me learn some stochastic calculus\ldots

%
%

\newpage
\section*{Solutions to exercises}

\begingroup
\newcommand{\separator}{
 \par
 \nointerlineskip
 \vskip5pt
 \hbox to\hsize{
  \hskip0.25\textwidth
  \leaders\hrule\hfil
  \hskip0.25\textwidth
 }
 \nointerlineskip
 \vskip10pt
}
\newcommand{\doubleseparator}{
 \par
 \nointerlineskip
 \vskip5pt
 \hbox to\hsize{
  \hskip0.35\textwidth
  \leaders\hrule height 1pt\hfil
  \hskip0.35\textwidth
 }
 \nointerlineskip
 \vskip10pt
}
\scriptsize

\newcounter{backupsection}
\setcounter{backupsection}{\value{section}}
\setcounter{section}{0}

{\stepcounter{section}

\doubleseparator

\begin{problem}
We have $P_y(x) = \lambda^d P_{\lambda y}(\lambda x)$. Thus, if $\lambda > 0$ and $f_\lambda(x) = f(\lambda x)$ and $u_\lambda(x, y)$ is the harmonic extension of $f_\lambda$, then
\[
 u_\lambda(x, y) = \int_{\R^d} P_y(x - x') f_\lambda(x') dx' = \int_{\R^d} \lambda^d P_{\lambda y}(\lambda x - \lambda x') f(\lambda x') dx' .
\]
Substituting $x'' = \lambda x'$, we obtain
\[
 u_\lambda(x, y) = \int_{\R^d} P_{\lambda y}(\lambda x - x'') f(x'') dx'' = u(\lambda x, \lambda y) .
\]
In particular,
\[
 K f_\lambda(x) = -\frac{\partial}{\partial y} u_\lambda(x, 0) = -\lambda \, \frac{\partial}{\partial y} u(\lambda x, 0) = \lambda K f(\lambda x) ,
\]
as long as all integrals and partial derivatives are well-defined.
\end{problem}

\separator

\begin{problem}\scriptsize
(a)~We have
\[
 \frac{\partial^2 P_y}{\partial x_k^2}(x) = -\frac{\partial}{\partial x_k} \, \frac{(d + 1) c_d x_k y}{(|x|^2 + y^2)^{(d + 3)/2}} = \frac{(d + 1) (d + 3) c_d x_k^2 y}{(|x|^2 + y^2)^{(d + 5)/2}} - \frac{(d + 1) c_d y}{(|x|^2 + y^2)^{(d + 3)/2}} \, ,
\]
and
\[
 \frac{\partial^2 P_y}{\partial y^2}(x) = \frac{\partial}{\partial y} \, \frac{c_d}{(|x|^2 + y^2)^{(d + 1)/2}} - \frac{\partial}{\partial y} \, \frac{(d + 1) c_d y^2}{(|x|^2 + y^2)^{(d + 3)/2}} = \frac{(d + 1) (d + 3) c_d y^3}{(|x|^2 + y^2)^{(d + 5)/2}} - \frac{3 (d + 1) c_d y}{(|x|^2 + y^2)^{(d + 3)/2}} \, .
\]
Therefore,
\[
 \Delta_{x,y} P_y(x) = \frac{(d + 1) (d + 3) c_d |x|^2 y}{(|x|^2 + y^2)^{(d + 5)/2}} - \frac{d (d + 1) c_d y}{(|x|^2 + y^2)^{(d + 3)/2}} + \frac{(d + 1) (d + 3) c_d y^3}{(|x|^2 + y^2)^{(d + 5)/2}} - \frac{3 (d + 1) c_d y}{(|x|^2 + y^2)^{(d + 3)/2}} \, .
\]
Rearranging the terms, we find that
\[
 \Delta_{x,y} P_y(x) = \frac{(d + 1) (d + 3) c_d (|x|^2 + y^2) y}{(|x|^2 + y^2)^{(d + 5)/2}} - \frac{(d + 1) (d + 3) c_d y}{(|x|^2 + y^2)^{(d + 3)/2}} = 0 .
\]

\noindent(b)~Note that $P_y(x)$ is an integrable function of $x \in \R^d$ for every $y > 0$. We set $f(x) = 1$ and $\lambda = 1 / y$ in the solution of the previous problem. We see that $u_\lambda(0, y) = u(0, 1)$. But $f_\lambda(x) = 1 = f(x)$, and so $u_\lambda(x, y) = u(x, y)$. Hence, $u(0, y) = u(0, 1)$. In other words, $\int_{\R^d} P_y(x) dx = \int_{\R^d} P_1(x) dx$. Choosing $c_d$ so that the latter integral is one, we find that $\int_{\R^d} P_y(x) dx = 1$ for every $y > 0$.

\noindent(c)~In order to prove that $u(x, y)$ is harmonic, we only need to observe that we can change the order of the integral and the application of the Laplace operator $\Delta_{x, y}$. Indeed: $P_y$ and all first and second order partial derivatives of $P_y$ are bounded by a constant times $P_y$, and so we may apply the dominated convergence theorem and change the order of partial derivatives and the integral. 

It remains to show that $u(x, y)$ converges to $f(x)$ as $y \to 0^+$. In fact, we have uniform convergence, as long as $f$ is bounded and uniformly continuous. This follows from a general result on approximate identities and the observation that the family of functions $P_y(\cdot)$ is an approximate identity as $y \to 0^+$; we omit the details.

\noindent(d)~We have $c_d = \pi^{-(d + 1)/2} \Gamma(\tfrac{d + 1}{2})$. There is a slick proof of this fact, which involves Gaussian kernels. By the definition of the Gamma function and substitution $s = (|x|^2 + y^2) / (4 t)$,
\[
 P_y(x) = \frac{1}{\pi^{(d + 1)/2}} \, \frac{y}{(|x|^2 + y^2)^{(d + 1) / 2}} \int_0^\infty s^{(d - 1) / 2} e^{-s} ds = \frac{y}{(4 \pi)^{(d + 1)/2}} \int_0^\infty t^{-(d + 3)/2} e^{-(|x|^2 + y^2) / (4 t)} dt .
\]
Integrating over $x \in \R^d$, substituting $r = y^2 / (4 t)$, and using $\Gamma(\tfrac12) = \sqrt{\pi}$, we obtain
\[
 \int_{\R^d} P_y(x) dx = \frac{y}{\sqrt{4 \pi}} \int_0^\infty t^{-3/2} e^{-y^2 / (4 t)} dt = \frac{1}{\sqrt{\pi}} \int_0^\infty r^{-1/2} e^{-r} dr = 1 .
\]
\end{problem}

\separator

\begin{problem}
We symmetrise the integral defining $u(x, y)$ and use dominated convergence theorem. Since $P_y(-x) = P_y(x)$, we have
\[
 u(x, y) = \int_{\R^d} f(x - z) P_y(z) dz = \int_{\R^d} f(x + z) P_y(z) dz .
\]
Hence,
\[
 u(x, y) = \frac{1}{2} \int_{\R^d} (f(x + z) + f(x - z)) P_y(z) dz .
\]
It follows that
\[
 K f(x) = \lim_{y \to 0^+} \frac{f(x) - u(x, y)}{y} = \lim_{y \to 0^+} \frac{1}{2} \int_{\R^d} (2 f(x) - f(x + z) - f(x - z)) \, \frac{P_y(z)}{y} \, dz .
\]
Since $f$ is smooth and bounded, we have $|2 f(x) - f(x + z) - f(x - z)| \le C \min\{1, |z|^2\}$ for some $C$. On the other hand, $0 \le P_y(z) / y \le c_d |z|^{-d - 1}$. Therefore, the dominated convergence theorem implies that
\[
 K f(x) = \frac{1}{2} \int_{\R^d} (2 f(x) - f(x + z) - f(x - z)) \biggl( \lim_{y \to 0^+} \frac{P_y(z)}{y} \biggr) dz ,
\]
which proves the first of the desired equalities. The other one follows by reversing the symmetrisation: if $B_\eps$ denotes the centred ball of radius $\eps$, then, by the same argument as in the first lines of the solution,
\[
 \frac{c_d}{2} \int_{\R^d \setminus B_\eps} \frac{2 f(x) - f(x + z) - f(x - z)}{|z|^{d + 1}} \, dz = c_d \int_{\R^d \setminus B_\eps} \frac{f(x) - f(x - z)}{|z|^{d + 1}} \, dz = c_d \int_{\R^d \setminus (x + B_\eps)} \frac{f(x) - f(x')}{|x - x'|^{d + 1}} \, dx' .
\]
We complete the proof by passing to the limit as $\eps \to 0^+$.
\end{problem}

\separator

\begin{problem}
The solution is very similar to that of exercise~1.2(d), but it requires the following observation. For every $a, b > 0$, substitution $u = a \sqrt{v} - b / \sqrt{v}$ yields
\[
 \sqrt{\pi} = \int_{-\infty}^\infty e^{-u^2} du = \frac{1}{2} \int_0^\infty (a v^{-1/2} + b v^{-3/2}) e^{-(a \sqrt{v} - b / \sqrt{v})^2} dv .
\]
Integrals of the two terms on the right-hand side are equal: if we set $v = b^2 / (a^2 w)$, then
\[
 \int_0^\infty a v^{-1/2} e^{-(a \sqrt{v} - b / \sqrt{v})^2} dv = \int_0^\infty b w^{-3/2} e^{-(b / \sqrt{w} - a \sqrt{w})^2} dw .
\]
Thus,
\[
 \sqrt{\pi} = \int_0^\infty b v^{-3/2} e^{-(a \sqrt{v} - b / \sqrt{v})^2} dv .
\]
Recall from the solution of exercise~1.2(d) that
\[
 P_y(x) = \frac{y}{(4 \pi)^{(d + 1)/2}} \int_0^\infty t^{-(d + 3)/2} e^{-(|x|^2 + y^2) / (4 t)} dt .
\]
Integrating over $x \in \R^d$, we obtain
\[
 \int_{\R^d} e^{i \xi x} P_y(x) dx = \frac{y}{\sqrt{4 \pi}} \int_0^\infty t^{-3/2} e^{-y^2 / (4 t)} e^{-t |\xi|^2} dt = \frac{e^{-y |\xi|}}{\sqrt{\pi}} \int_0^\infty \frac{y}{2 t^{3/2}} e^{-(|\xi| \sqrt{t} - y / \sqrt{4 t})^2} dt = e^{-y |\xi|} .
\]
\end{problem}

\separator

\begin{problem}
(a)~The process $Y^{(\lambda)}_t = \lambda^{-1/2} Y_{\lambda t}$ has the same law as the process $Y_t$, and so its local time has the same law as $L_t$. On the other hand, the corresponding approximation to the local time $L_t^{(\lambda, \delta)}$ is equal to $\lambda^{-1} L_{\lambda t}^{(\lambda^{1/2} \delta)}$, and since the local time of $Y^{(\lambda)}_t$ is the limit of $(2 \delta)^{-1} L_t^{(\lambda, \delta)}$, it is equal to $\lambda^{-1/2} L_{\lambda t}$. It follows that $\lambda^{-1/2} L_{\lambda t}$ has the same law as $L_t$.

\noindent(b)~This follows directly from the definition: $T_{\lambda s} = \max\{t \ge 0 : \lambda^{-1} L_t \le s\}$ is equal in law to $\max\{t \ge 0 : L_{\lambda^{-2} t} \le s\} = \lambda^2 T_s$.

\noindent(c)~We know that $T_{\lambda s}$ is equal in law to $\lambda^2 T_s$, and $X_{\lambda^2 t}$ is equal in law to $\lambda X_t$. Thus, the processes $X(T_{\lambda s})$, $X(\lambda^2 T_s)$ and $\lambda X(T_s)$ are all equal in law, as desired. Symmetry is clear.
\end{problem}
}

{\stepcounter{section}

\doubleseparator

\begin{problem}
(a)~The calculation is not much more complicated than in the classical case $\alpha = 1$. We have
\[
 \frac{\partial^2 P_y}{\partial x_k^2}(x) = -\frac{\partial}{\partial x_k} \, \frac{(d + \alpha) c_{d, \alpha} x_k y^\alpha}{(|x|^2 + y^2)^{(d + \alpha + 2)/2}} = \frac{(d + \alpha) (d + \alpha + 2) c_{d, \alpha} x_k^2 y^\alpha}{(|x|^2 + y^2)^{(d + \alpha + 4)/2}} - \frac{(d + \alpha) c_{d, \alpha} y^\alpha}{(|x|^2 + y^2)^{(d + \alpha + 2)/2}} \, ,
\]
\[
 \frac{\partial P_y}{\partial y}(x)  = \frac{\alpha c_{d, \alpha} y^{\alpha - 1}}{(|x|^2 + y^2)^{(d + \alpha)/2}} - \frac{(d + \alpha) c_{d, \alpha} y^{\alpha + 1}}{(|x|^2 + y^2)^{(d + \alpha + 2)/2}} \, ,
\]
and
\[
 \frac{\partial^2 P_y}{\partial y^2}(x) = \frac{(d + \alpha) (d + \alpha + 2) c_{d, \alpha} y^{\alpha + 2}}{(|x|^2 + y^2)^{(d + \alpha + 4)/2}} - \frac{(2 \alpha + 1) (d + \alpha) c_{d, \alpha} y^\alpha}{(|x|^2 + y^2)^{(d + \alpha + 2)/2}} + \frac{\alpha (\alpha - 1) c_{d, \alpha} y^{\alpha - 2}}{(|x|^2 + y^2)^{(d + \alpha)/2}} \, .
\]
Therefore, after simple rearrangement of terms,
\[
 \Delta_{x,y} P_y(x) = \frac{(d + \alpha) (d + \alpha + 2) c_{d, \alpha} (|x|^2 + y^2) y^\alpha}{(|x|^2 + y^2)^{(d + \alpha + 4)/2}} - \frac{(d + \alpha) (d + 2 \alpha + 1) c_{d, \alpha} y^\alpha}{(|x|^2 + y^2)^{(d + \alpha + 2)/2}} + \frac{\alpha (\alpha - 1) c_{d, \alpha} y^{\alpha - 2}}{(|x|^2 + y^2)^{(d + \alpha)/2}} \, .
\]
Simplification leads to
\[
 \Delta_{x,y} P_y(x) = -\frac{(\alpha - 1) (d + \alpha) c_{d, \alpha} y^\alpha}{(|x|^2 + y^2)^{(d + \alpha + 2)/2}} + \frac{\alpha (\alpha - 1) c_{d, \alpha} y^{\alpha - 2}}{(|x|^2 + y^2)^{(d + \alpha)/2}} = \frac{\alpha - 1}{y} \, \frac{\partial P_y}{\partial y}(x) .
\]

\noindent(b)~We can argue as in the classical case $\alpha = 1$, or observe that $P_y(x) = y^{-d} P_1(x / y)$, and so $\int_{\R^d} P_y(x) dx = \int_{\R^d} P_1(x) dx$. Choosing $c_{d, \alpha}$ so that the latter integral is one, we find that $\int_{\R^d} P_y(x) dx = 1$ for every $y > 0$.

\noindent(c)~We proceed as in the classical case. By the definition of the Gamma function and substitution $s = (|x|^2 + y^2) / (4 t)$,
\[
 P_y(x) = \frac{c_{d, \alpha} y^\alpha}{\Gamma(\tfrac{d + \alpha}{2}) (|x|^2 + y^2)^{(d + \alpha) / 2}} \int_0^\infty s^{(d + \alpha) / 2 - 1} e^{-s} ds = \frac{c_{d, \alpha} y^\alpha}{2^{d + \alpha} \Gamma(\tfrac{d + \alpha}{2})} \int_0^\infty t^{-(d + \alpha) / 2 - 1} e^{-(|x|^2 + y^2) / (4 t)} dt .
\]
Integrating over $x \in \R^d$, substituting $r = y^2 / (4 t)$, and once again using the definition of the Gamma function, we obtain
\[
 \int_{\R^d} P_y(x) dx = \frac{\pi^{d / 2} c_{d, \alpha} y^\alpha}{2^\alpha \Gamma(\tfrac{d + \alpha}{2})} \int_0^\infty t^{-\alpha/2 - 1} e^{-y^2 / (4 t)} dt = \frac{\pi^{d / 2} c_{d, \alpha}}{\Gamma(\tfrac{d + \alpha}{2})} \int_0^\infty r^{\alpha / 2 - 1} e^{-r} dr = \frac{\pi^{d / 2} \Gamma(\tfrac{\alpha}{2}) c_{d, \alpha}}{\Gamma(\tfrac{d + \alpha}{2})} .
\]
Thus, we have $c_{d,\alpha} = \pi^{-d / 2} \Gamma(\tfrac{d + \alpha}{2}) / \Gamma(\tfrac{\alpha}{2})$. 
\end{problem}

\separator

\begin{problem}
The derivation is exactly the same as in the classical case $\alpha = 1$, except that the estimate $0 \le P_y(z) / y \le c_d |z|^{-d - 1}$ now becomes $0 \le P_y(z) / y^\alpha \le c_{d, \alpha} |z|^{-d - \alpha}$. This does not affect the applicability of the dominated convergence theorem, and so the Dirichlet-to-Neumann operator, defined by the formula
\[
 K f(x) = \lim_{y \to 0^+} \frac{f(x) - u(x, y)}{y^\alpha} = \lim_{y \to 0^+} \frac{1}{2} \int_{\R^d} (2 f(x) - f(x + z) - f(x - z)) \, \frac{P_y(z)}{y^\alpha} \, dz ,
\]
is given by
\[
 K f(x) = \frac{1}{2} \int_{\R^d} (2 f(x) - f(x + z) - f(x - z)) \biggl( \lim_{y \to 0^+} \frac{P_y(z)}{y^\alpha} \biggr) dz = \frac{c_{d, \alpha}}{2} \int_{\R^d} \frac{2 f(x) - f(x + z) - f(x - z)}{|z|^{d + \alpha}} \, dz .
\]
In the same way one then shows that
\[
 K f(x) = \lim_{\eps \to 0^+} c_{d, \alpha} \int_{\R^d \setminus (x + B_\eps)} \frac{f(x) - f(x')}{|x - x'|^{d + \alpha}} \, dx' .
\]
This is \emph{nearly} the same expression as the one for $(-\Delta)^{\alpha / 2} f(x)$: the constant $c_{d, \alpha} = \pi^{-d / 2} \Gamma(\tfrac{d + \alpha}{2}) / \Gamma(\tfrac{\alpha}{2})$ is wrong! The definition of $(-\Delta)^{\alpha / 2}$ has $c_{d, \alpha} = -2^\alpha \pi^{-d / 2} \Gamma(\tfrac{d + \alpha}{2}) / \Gamma(-\tfrac{\alpha}{2})$, so that
\[
 K f(x) = -\frac{\Gamma(-\tfrac{\alpha}{2})}{2^\alpha \Gamma(\tfrac{\alpha}{2})} \, (-\Delta)^{\alpha/2} f(x) .
\]
It is no coincidence that the constant on the right-hand side does not depend on $d$.
\end{problem}

\separator

\begin{problem}
If $v(x, z) = u(x, z^{1 / \alpha})$, then
\[
 \frac{\partial v}{\partial z}(x, z) = \alpha^{-1} z^{1 / \alpha - 1} \, \frac{\partial u}{\partial y}(x, z^{1 / \alpha})
\]
and
\[
 \frac{\partial^2 v}{\partial z^2}(x, z) = \alpha^{-2} z^{2 / \alpha - 2} \, \frac{\partial^2 u}{\partial y^2}(x, z^{1 / \alpha}) + \alpha^{-1} (\alpha^{-1} - 1) z^{1 / \alpha - 2} \, \frac{\partial u}{\partial y}(x, z^{1 / \alpha}) .
\]
Obviously, $\Delta_x v(x, z) = \Delta_x u(x, z^{1 / \alpha})$. Therefore,
\[
 \alpha^{-2} z^{2 / \alpha - 2} \Delta_x v(x, z) + \frac{\partial^2 v}{\partial z^2} (x, z) = \alpha^{-2} z^{2 / \alpha - 2} \biggl(\Delta_x u(x, z^{1 / \alpha}) + \frac{\partial^2 u}{\partial y^2}(x, z^{1 / \alpha}) + \frac{1 - \alpha}{z^{1 / \alpha}} \, \frac{\partial u}{\partial y}(x, z^{1 / \alpha}) \biggr) = 0 .
\]
\end{problem}

\separator

\begin{problem}
This problem requires some knowledge of the Bessel $K$ function, which we assume here. We have
\[
 \fourier f(\xi) = \fourier u(\xi, 0) = \lim_{y \to 0^+} \fourier u(\xi, y) = c_\alpha \fourier f(\xi) \lim_{y \to 0^+} (y |\xi|)^{\alpha/2} K_{\alpha/2} (y |\xi|) ,
\]
and so
\[
 1 = c_\alpha \lim_{s \to 0^+} s^{\alpha/2} K_{\alpha/2}(s) = c_\alpha 2^{\alpha/2 - 1} \Gamma(\tfrac{\alpha}{2}) .
\]
Therefore, $c_\alpha = 2^{1 - \alpha/2} / \Gamma(\tfrac{\alpha}{2})$.

It is perhaps worthwhile to note that in a similar way we can recover the constant in the identity
\[
 K f(x) = -\frac{\Gamma(-\tfrac{\alpha}{2})}{2^\alpha \Gamma(\tfrac{\alpha}{2})} \, (-\Delta)^{\alpha/2} f(x) .
\]
Indeed, with $\nu = \tfrac{\alpha}{2} \in (0, 1)$, we have
\[
 s^\nu K_\nu(s) = 2^{\nu - 1} \Gamma(\nu) + 2^{-\nu - 1} \Gamma(-\nu) s^{2 \nu} + o(s^{2 \nu})
\]
as $s \to 0^+$. Since $c_\alpha = (2^{\nu - 1} \Gamma(\nu))^{-1}$, we obtain
\[
 \fourier K f(\xi) = |\xi|^\alpha \fourier f(\xi) \lim_{s \to 0^+} \frac{1 - c_\alpha s^\nu K_\nu(s)}{s^{2 \nu}} = -\frac{\Gamma(-\nu)}{2^{2 \nu} \Gamma(\nu)} \, |\xi|^\alpha \fourier f(\xi) ,
\]
as desired. Note that with this approach, it is clear that the constant does not depend on the dimension.
\end{problem}

\separator

\begin{problem}
(a)~Suppose that $Y_t$ is the solution of $dY_t = dW_t + (1 - \alpha) dY_t / Y_t$, where $W_t$ is the Brownian motion. The process $Y^{(\lambda)}_t = \lambda^{-1/2} Y_{\lambda t}$ satisfies the stochastic differential equation
\[
 dY^{(\lambda)}_t = \lambda^{-1/2} d[W_{\lambda t}] + \frac{1 - \alpha}{\sqrt{\lambda} Y_t} \, d[\lambda t] = dW_t^{(\lambda)} + \frac{1 - \alpha}{Y^{(\lambda)}_t} \, dt ,
\]
where $W_t^{(\lambda)} = \lambda^{-1/2} W_{\lambda t}$ is the Brownian motion. It follows that $Y^{(\lambda)}_t$ is indeed the same Bessel process, as desired.

The above proof is not entirely correct when $\alpha > 1$, because then $Y_t$ is too singular near $0$. A fully rigorous argument uses the squared Bessel process $Q_t = Y_t^2$, which solves
\[
 dQ_t = 2 \sqrt{Q_t} \, dW_t + (2 - \alpha) dt .
\]
For $\lambda > 0$, $Q_t^{(\lambda)} = \lambda^{-1} Q_{\lambda t}$ satisfies the same equation driven by $W_t^{(\lambda)} = \lambda^{-1/2} W_{\lambda t}$. Uniqueness in law proves that $Q_t^{(\lambda)}$ and $Q_t$ are equal in law, and taking square roots proves the claim.

\noindent(b)~The argument is very similar to the solution of exercise~1.5(a). The approximation $L_t^{(\lambda, \delta)}$ to the local time of $Y^{(\lambda)}_t$ is equal to $\lambda^{-1} L_{\lambda t}^{(\lambda^{1/2} \delta)}$, and the local time of $Y_t^{(\lambda)}$ is the limit of $\tfrac12 \delta^{-2 + \alpha} L_t^{(\lambda, \delta)}$. It follows that the local time of $Y_t^{(\lambda)}$ is equal to $\lambda^{-\alpha/2} L_{\lambda t}$. It follows that $\lambda^{-\alpha/2} L_{\lambda t}$ has the same law as $L_t$.

\noindent(c)~As in exercise~1.5(b), by definition, $T_{\lambda s} = \max\{t \ge 0 : \lambda^{-1} L_t \le s\}$ is equal in law to $\max\{t \ge 0 : L_{\lambda^{-2 / \alpha} t} \le s\} = \lambda^{2 / \alpha} T_s$.

\noindent(d)~We follow the solution of exercise~1.5(c). We know that $T_{\lambda s}$ is equal in law to $\lambda^{2 / \alpha} T_s$, and $X_{\lambda^{2 / \alpha} t}$ is equal in law to $\lambda^{1 / \alpha} X_t$. Thus, the processes $X(T_{\lambda s})$, $X(\lambda^{2 / \alpha} T_s)$ and $\lambda^{1 / \alpha} X(T_s)$ are all equal in law, as desired. Symmetry is clear.
\end{problem}
}

{\stepcounter{section}

\doubleseparator

\begin{problem}
We need to find a positive, nonincreasing solution $\ph = \ph_\lambda$ of the ODE
\[
 \ph''(y) = \lambda a(y) \ph(y) = \lambda \ind_{(0, 1)}(y) \ph(y) ,
\]
with $\ph(0) = 1$. Clearly, $\ph$ is linear, positive and nonincreasing on $[1, \infty)$, and hence constant on $[1, \infty)$. It follows that $\ph'(1) = 0$, and so $\ph(y) = c \cosh(\sqrt{\lambda} (1 - y))$ for $y \in [0, 1]$. Since $\ph(0) = 1$, we have $c = (\cosh \sqrt{\lambda})^{-1}$. Now,
\[
 \mu(\lambda) = -\ph'(0) = c \sqrt{\lambda} \sinh \sqrt{\lambda} = \sqrt{\lambda} \tanh \sqrt{\lambda} \, .
\]
\end{problem}

\separator

\begin{problem}
Similarly, we find a solution $\ph = \ph_\lambda$ on $(0, 1)$ of the ODE
\[
 \ph''(y) = \lambda a(y) \ph(y) = \lambda \ph(y) ,
\]
with $\ph(0) = 1$ and $\ph(1) = 0$. For $\lambda > 0$, $\ph(y) = c \sinh(\sqrt{\lambda} (1 - y))$ for $y \in [0, 1]$, and since $\ph(0) = 1$, we have $c = (\sinh \sqrt{\lambda})^{-1}$. We find that
\[
 \mu(\lambda) = -\ph'(0) = c \sqrt{\lambda} \cosh \sqrt{\lambda} = \sqrt{\lambda} \coth \sqrt{\lambda} \, .
\]
When $\lambda = 0$, we have $\ph(y) = 1 - y$ and $\mu(0) = 1$.
\end{problem}

\separator

\begin{problem}
This is quite straightforward: a positive, nonincreasing solution $\ph = \ph_\lambda$ of the ODE
\[
 \ph''(y) = \lambda a(y) \ph(y) = 0
\]
is necessarily constant, and since $\ph(0) = 1$, we have $\ph(y) = 1$ and
\[
 \mu(\lambda) = -\ph'(0) = 0 .
\]
\end{problem}

\separator

\begin{problem}
This is very similar: a solution $\ph = \ph_\lambda$ of the ODE
\[
 \ph''(y) = \lambda a(y) \ph(y) = 0
\]
is necessarily linear, and since $\ph(0) = 1$ and $\ph(1) = 0$, we have $\ph(y) = 1 - y$ and
\[
 \mu(\lambda) = -\ph'(0) = 1 .
\]
\end{problem}

\separator

\begin{problem}
We find a positive, nonincreasing distributional solution $\ph = \ph_\lambda$ of the ODE
\[
 \ph''(y) = \lambda a(y) \ph(y) ,
\]
where $a(y)$ is the Dirac delta at $1$. Clearly, $\ph$ is linear on $[0, 1]$, constant on $[1, \infty)$. Furthermore, $\ph''$ has an atom at $1$ with mass $\ph'(1^+) - \ph'(1^-)$, and hence $\ph'(1^+) - \ph'(1^-) = \lambda \ph(1)$. Clearly, if we write $\mu$ for $\mu(\lambda)$, then $\ph(0) = 1$ and $\ph'(0) = -\mu$. It follows that $\ph(y) = 1 - \mu y$ for $y \in [0, 1]$, and $\ph(y) = 1 - \mu$ for $y \in [1, \infty)$. We get
\[
 \mu = 0 - (-\mu) = \ph'(1^+) - \ph'(1^-) = \lambda \ph(1) = \lambda (1 - \mu) ,
\]
which leads to $\mu(\lambda) = \mu = \lambda / (\lambda + 1)$.
\end{problem}

\separator

\begin{problem}
This is no different than the solution of exercise~3.3: our ODE is simply $\ph''(y) = 0$, and hence $\ph(y) = 1$. But according to our generalised definition of $\ph'(0)$, we have
\[
 \mu(\lambda) = -\ph'(0) = -(\ph'(0^+) - \lambda a(\{0\}) \ph(0)) = -(0 - \lambda \cdot 1 \cdot 1) = \lambda . 
\]
\end{problem}

\separator

\begin{problem}
Here we need to find a positive nonincreasing solution $\ph = \ph_\lambda$ of the ODE
\[
 \ph''(y) = \lambda a(y) \ph(y) = \frac{\lambda \ph(y)}{(1 + 2 y)^2} \, ,
\]
with $\ph(0) = 1$. This is not as straightforward as the previous problems. Substituting $1 + 2 y = s$, $\psi(s) = \ph(y) = \ph(\tfrac12 (s - 1))$, we obtain
\[
 \psi''(s) = \frac{\ph''(y)}{4} = \frac{\lambda \ph(y)}{4 (1 + 2 y)^2} = \frac{\lambda \psi(s)}{4 s^2} \, ,
\]
with $\psi$ nonincreasing and $\psi(1) = 1$. It is therefore natural to search for solutions of the form $\psi(s) = s^p$, with $p \le 0$. With this choice,
\[
 p (p - 1) s^{p - 2} = \psi''(s) = \frac{\lambda \psi(s)}{4 s^2} = \frac{\lambda}{4} \, s^{p - 2} ,
\]
and so $4 p (p - 1) = \lambda$, that is, $(2 p - 1)^2 = \lambda + 1$, or $p = \tfrac12 (1 \pm \sqrt{\lambda + 1})$. Since $p \le 0$, we choose the minus sign, and it follows that
\[
 \ph(y) = (1 + 2 y)^p = (1 + 2 y)^{(1 - \sqrt{\lambda + 1}) / 2} .
\]
We conclude that
\[
 \mu(\lambda) = -\ph'(0) = \sqrt{\lambda + 1} - 1 .
\]
\end{problem}

\separator

\begin{problem}
This is very similar to the previous problem: we are looking for a positive nonincreasing solution $\ph = \ph_\lambda$ of the ODE
\[
 \ph''(y) = \lambda a(y) \ph(y) = \frac{\lambda \ph(y)}{(1 - 2 y)^2}
\]
on $(0, \tfrac12)$, with $\ph(0) = 1$. Substituting $1 - 2 y = s$, $\psi(s) = \ph(y) = \ph(\tfrac12 (1 - s))$, we obtain the same ODE as before,
\[
 \psi''(s) = \frac{\ph''(y)}{4} = \frac{\lambda \ph(y)}{4 (1 - 2 y)^2} = \frac{\lambda \psi(s)}{4 s^2} \, ,
\]
but with $\psi$ nondecreasing and $\psi(1) = 1$. We have the same set of solutions, but this time we choose the plus sign. That is, we have $\psi(s) = s^p$ with $p = \tfrac12 (1 + \sqrt{\lambda + 1})$. It follows that
\[
 \ph(y) = (1 - 2 y)^p = (1 - 2 y)^{(1 + \sqrt{\lambda + 1}) / 2} ,
\]
and therefore
\[
 \mu(\lambda) = -\ph'(0) = \sqrt{\lambda + 1} + 1 .
\]
The other solution of our ODE, with $p = \tfrac12 (1 - \sqrt{\lambda + 1})$, is indeed unbounded near $y = \tfrac12$ (when $\lambda > 0$), and so $\ph$ described above is the only bounded solution.
\end{problem}

\separator

\begin{problem}
Once we know the right operator, the proof is fairly easy, but it is more interesting to understand how one can actually derive the solution.

\emph{Part I. Informal derivation.} Let us start with something that we know very well: the operator
\[
 L = \Delta_x + \frac{\partial^2}{\partial y^2}
\]
leads to the Dirichlet-to-Neumann map
\[
 K = \sqrt{-\Delta_x} .
\]
It is fairly easy to believe that $-\Delta_x$ can be replaced by an abstract positive self-adjoint operator on a Hilbert space: we just need to refer to the spectral representation of our operator instead of using the Fourier transform. And it is far easier to see that we can replace $-\Delta_x$ by $-\Delta_x + 1$: the proof is exactly the same! That is, the operator
\[
 L = \Delta_x - 1 + \frac{\partial^2}{\partial y^2}
\]
leads to the desired Dirichlet-to-Neumann map
\[
 K = \sqrt{-\Delta_x + 1} .
\]
But this $L$ is not of the desired form! It is, however, quite standard to transform it into an operator that we are looking for without changing the Dirichlet-to-Neumann map. First, we apply what is sometimes known as the Doob transform: we find a function $\psi(y) = \cosh y$ which satisfies $(-1 + \frac{\partial^2}{\partial y^2}) \psi(y) = 0$ and $\psi'(0) = 0$, and we observe that
\[
 L [u(x, y) \psi(y)] = \Delta_x u(x, y) \psi(y) - u(x, y) \psi(y) + \frac{\partial^2 u}{\partial y^2}(x, y) \psi(y) + 2 \, \frac{\partial u}{\partial y}(x, y) \, \frac{\partial \psi}{\partial y}(y) + u(x, y) \, \frac{\partial^2 \psi}{\partial y^2}(y) .
\]
After simplification, we arrive at
\[
 (\psi(y))^{-1} L [u(x, y) \psi(y)] = \Delta_x u(x, y) + \frac{\partial^2 u}{\partial y^2}(x, y) + 2 \tanh y \, \frac{\partial u}{\partial y}(x, y) .
\]
Next, we substitute $s = \tanh y$ and $u(x, y) = v(x, s) = v(x, \tanh y)$, so that
\[
 (\psi(y))^{-1} L [u(x, y) \psi(y)] = \Delta_x v(x, s) + \frac{\partial^2 v}{\partial s^2}(x, s) \, \frac{1}{(\cosh y)^4} + \frac{\partial v}{\partial s}(x, s) \, \frac{-2 \sinh y}{(\cosh y)^3} + 2 \tanh y \, \frac{\partial v}{\partial s}(x, s) \, \frac{1}{(\cosh y)^2} \, .
\]
Simplifying the right-hand side leads to
\[
 (\psi(y))^{-1} L [u(x, y) \psi(y)] = \Delta_x v(x, s) + (1 - s^2)^2 \, \frac{\partial^2 v}{\partial s^2}(x, s) .
\]
It remains to divide both sides by $(1 - s^2)^2$, and we get the final answer: the operator we are looking for maps $v(x, s)$ to
\[
 \frac{1}{(1 - s^2)^2} \, (\psi(y))^{-1} L [u(x, y) \psi(y)] = \frac{1}{(1 - s^2)^2} \, \Delta_x v(x, s) + \frac{\partial^2 v}{\partial s^2}(x, s) .
\]
And only now the actual solution begins. We could use $y$ instead of $s$, but just to be consistent with the above derivation, we will stick to $s$.

\emph{Part 2. Rigorous solution.}
We claim that the operator
\[
 L = \frac{1}{(1 - s^2)^2} \, \Delta_x + \frac{\partial^2}{\partial s^2}
\]
on $\R^d \times (0, 1)$, with zero boundary data at $s = 1$, corresponds to the desired Dirichlet-to-Neumann operator $K = \sqrt{-\Delta + 1}$. In order to prove that, we need to find a nonincreasing positive solution $\ph = \ph_\lambda$ of the ODE
\[
 \ph''(s) = \lambda (1 - s^2)^{-2} \ph(s) ,
\]
which satisfies $\ph(0) = 1$. Perhaps the easiest way is to reverse the steps that have led us to the above ODE.

First, we substitute $s = \tanh y$ and $\psi(y) = \ph(s) = \ph(\tanh y)$. The above equation transforms to
\[
 \psi''(y) = \frac{\ph''(s)}{(\cosh y)^4} - \frac{2 \sinh y \, \ph'(s)}{(\cosh y)^3} = (1 - s^2)^2 \ph''(s) - \frac{2 \sinh y \, \ph'(s)}{(\cosh y)^3} = \lambda \ph(s) - \frac{2 \sinh y \, \ph'(s)}{(\cosh y)^3} \, .
\]
Since $\psi'(y) = \ph'(s) (\cosh y)^{-2}$, we find that
\[
 \psi''(y) = \lambda \psi(y) - 2 \tanh y \, \psi'(y) .
\]
Next, if we denote $\eta(y) = \cosh y \, \psi(y)$, then
\[
 \eta''(y) = \cosh y \, (\psi(y) + 2 \tanh y \, \psi'(y) + \psi''(y)) = (\lambda + 1) \cosh y \, \psi(y) = (\lambda + 1) \eta(y) .
\]
Furthermore, $1 = \ph(0) = \psi(0) = \eta(0)$, and if $\eta$ is nonincreasing, then so is $\psi$ and $\ph$. Thus,
\[
 \eta(y) = e^{-\sqrt{\lambda + 1} \, y} , \qquad \psi(y) = \frac{e^{-\sqrt{\lambda + 1} \, y}}{\cosh y} \, , \qquad \ph(s) = \sqrt{1 - s^2} e^{-\sqrt{\lambda + 1} \, \artanh s} .
\]
In fact, we do not need the rather ugly expressions for $\psi$ and $\ph$: by the formula for $\eta$, we have
\[
 \mu(\lambda) = -\ph'(0) = -\psi'(0) = -\eta'(0) = \sqrt{\lambda + 1} ,
\]
as desired.
\end{problem}
}

{\stepcounter{section}

\doubleseparator

\begin{problem}
While there can be no `solution' to this problem, it is clear that the key property of $\ph$ that is needed here is that $\Xi_t \ph(Y_t)$ is a martingale until $Y_t$ reaches zero. Thus, one could approach the problem from the opposite direction: evaluate $d(\Xi_t \ph(Y_t))$ using all the stochastic machinery (without specifying $\ph$), and figure out that this is a stochastic differential of a martingale (up to $T_0$) only if $\ph$ satisfies Krein's ODE.
\end{problem}

\separator

\begin{problem}
(a)~The calculation is exactly the same, except that all integrals with respect to the local time $L_t$ and with respect to variable $r$ are gone, and the initial value is $\Xi_0 \Phi_0 = e^{-i \xi x} \ph(y)$ rather than $1$. That is,
\[
 \ex_Y[\Xi(T_0)] = \ex_Y[\Xi(T_0) \ph(Y(T_0))] = e^{-i \xi x} \ph(y) .
\]
Here, of course, $\ph(y) = \ph(|\xi|^2, y)$.

\noindent(b)~It suffices to observe that $\ex \exp(-i \xi X(T_0)) = \ex_Y[\ex_W[\exp(-i \xi X(T_0))]] = \ex_Y[\Xi(T_0)]$ (recall that $T_0$ and $W$ are independent) and use the result of item~(a).

\noindent(c)~We have seen that the Poisson kernel $P_y(x)$ is the distribution of a random variable, and hence it is nonnegative. Recall that $u(x, y)$ (as a function of $x$) is the convolution of the boundary values $f(x)$ and the Poisson kernel $P_y(x)$. Hence, if $f$ is nonnegative, then $u(x, y)$ is nonnegative, too.

\noindent(d)~The argument is very similar: the convolution of a compactly supported continuous function $f$ with a subprobability measure is continuous, with the same modulus of continuity as $f$. This gives equicontinuity of $u(\cdot, y)$ with respect to $y > 0$. Recall that $u(\cdot, y)$ is continuous with respect to $y$ as a function with values in $L^2$. This is sufficient to conclude that $u$ is indeed continuous. Indeed: suppose that $u$ has a discontinuity at $(x, y)$, that is, for some $\eps > 0$ we can have $|u(x, y) - u(x', y')| \ge 3 \eps$ with $(x', y')$ arbitrarily close to $(x, y)$. Choose $\delta > 0$ according to the definition of equicontinuity, and consider $(x', y')$ as above, with $|x - x'| < \eta$ and $|y - y'| < \eta$, where $\eta < \delta$. By equicontinuity, we have
\[
 |u(x'', y) - u(x'', y')| \ge |u(x, y) - u(x', y')| - |u(x'', y) - u(x, y)| - |u(x'', y') - u(x', y')| \ge 3 \eps - \eps - \eps = \eps
\]
if $|x'' - x| < \delta$ and $|x'' - x'| < \delta$. This condition is satisfied when $|x'' - \tfrac{1}{2} (x + x')| < \delta - \tfrac{1}{2} \eta$, and so in a ball of radius at least $\tfrac{1}{2} \delta$. Thus, the $L^2$ distance of $u(\cdot, y)$ and $u(\cdot, y')$ is greater than a constant times $\delta^{d / 2}$, a contradiction. 

(Similarly, if $f$ is additionally smooth, then $u(x, y)$ is smooth with respect to $x$; but, unfortunately, $u(x, y)$ may fail to be smooth with respect to $y$.)
\end{problem}

\separator

\begin{problem}
This is an open-ended question!

Suppose that reflection at the boundary is realised by forcing an upward jump when $Y_n = 0$. Let
\[
 h(\xi) = \frac{1}{d} \sum_{j = 1}^d \cos \xi_j
\]
be the characteristic function of a horizontal step, conditional on making one. Between consecutive vertical jumps in the interior, the number of horizontal steps is geometric on $\{0, 1, \ldots\}$ with success probability $\frac{1}{d + 1}$. Thus, the horizontal displacement in such a block has characteristic function
\[
 \ph(\xi) = \sum_{n = 0}^\infty \frac{1}{d + 1} \biggl(\frac{d}{d+1} \, h(\xi)\biggr)^n = \frac{1}{d + 1 - d h(\xi)} \, .
\]
The forced upward step from the boundary has no horizontal displacement and must be treated separately.

Let $\tau_0$ be the first hitting time of height $0$, and define
\[
 f_j(\xi) = \ex_{(0, j)} \exp(-i \xi X(\tau_0)) .
\]
We have $f_0(\xi) = 1$, and conditioning on the next vertical jump gives
\[
 f_j(\xi) = \tfrac12 \ph(\xi) (f_{j + 1}(\xi) + f_{j - 1}(\xi))
\]
for $j = 1, 2, \ldots$\, The bounded solution is
\[
 f_j(\xi) = \biggl(\frac{1 - \sqrt{1 - (\ph(\xi))^2}}{\ph(\xi)}\biggr)^j = \biggl(\frac{\ph(\xi)}{1 + \sqrt{1 - (\ph(\xi))^2}}\biggr)^j .
\]
By the strong Markov property at successive boundary visits, $Z_k$ is a random walk with independent identically distributed increments having characteristic function $f_1(\xi)$.

Note that as $\xi \to 0$, we have
\[
 h(\xi) = 1 - \frac{|\xi|^2}{2 d} + O(|\xi|^4) ,
\]
and so
\[
 f_1(\xi) = 1 - |\xi| + O(|\xi|^2)
\]
near zero, displaying the same scaling of the trace at large scales as for the Cauchy process.
\end{problem}
}

\setcounter{section}{\value{backupsection}}

\endgroup

%
%

%
%


\newpage
\begin{thebibliography}{00}

\bibitem{assing-herman}
\textsc{Sigurd Assing, John Herman,}
\emph{Extension technique for functions of diffusion operators: a stochastic approach,}
Electron. J.~Probab. 26 (2021), paper no. 67: 1--32.
\doi{10.1214/21-EJP624}

\bibitem{caffarelli-silvestre}
\textsc{Luis Caffarelli, Luis Silvestre,}
\emph{An extension problem related to the fractional Laplacian,}
Comm. Partial Differential Equations 32(7) (2007): 1245--1260.
\doi{10.1080/03605300600987306}

\bibitem{eckhardt-kostenko}
\textsc{Jonathan Eckhardt, Aleksey Kostenko,}
\emph{The inverse spectral problem for indefinite strings,}
Invent. Math. 204 (2016): 939--977.
\doi{10.1007/s00222-015-0629-1}

\bibitem{hauer-lee}
\textsc{Daniel Hauer, David Lee,}
\emph{Functional Calculus via the extension technique: a first hitting time approach,}
preprint (2021).
\arxiv{2101.11305}

\bibitem{debranges-1}
\textsc{Louis de Branges,}
\emph{Some Hilbert spaces of entire functions.}
Trans. Amer. Math. Soc. 96 (1960): 259--295.
\doi{10.1090/S0002-9947-1960-0133455-X}

\bibitem{debranges-2}
\textsc{Louis de Branges,}
\emph{Some Hilbert spaces of entire functions II.}
Trans. Amer. Math. Soc. 99 (1961): 118--152.
\doi{10.1090/S0002-9947-1961-0133456-2}

\bibitem{debranges-3}
\textsc{Louis de Branges,}
\emph{Some Hilbert spaces of entire functions III.}
Trans. Amer. Math. Soc. 100 (1961): 73--115.
\doi{10.1090/S0002-9947-1961-0133457-4}

\bibitem{debranges-4}
\textsc{Louis de Branges,}
\emph{Some Hilbert spaces of entire functions IV.}
Trans. Amer. Math. Soc. 105 (1962): 43--83.
\doi{10.1090/S0002-9947-1962-0143016-6}

\bibitem{debranges}
\textsc{Louis de Branges,}
\emph{Hilbert spaces of entire functions.}
Prentice-Hall Inc, Englewood Cliffs, 1968.

\bibitem{ito-mckean}
\textsc{Kiyosi Itô, Henry P.~McKean,}
\emph{Diffusion Processes and their Sample Paths.}
Grundlehren der mathematischen Wissenschaften~125,
Springer-Verlag, Berlin-Heidelberg-New York, 1974.

\bibitem{kotani-watanabe}
\textsc{Shinichi Kotani, Shinzo Watanabe,}
\emph{Krein's spectral theory of strings and generalized diffusion processes.}
In: M.~Fukushima (eds), \emph{Functional Analysis in Markov Processes,}
Lecture Notes in Mathematics 923,
Springer, Berlin, Heidelberg, 1982.
\doi{10.1007/BFb0093046}

\bibitem{krein}
\textsc{Mark G.~Krein,}
\emph{On a generalization of an investigation of Stieltjes,}
Dokl. Akad. Nauk SSSR 87 (1952): 881--884 (in Russian).

\bibitem{knight}
\textsc{Frank B.~Knight,}
\emph{Characterization of the Levy measures of inverse local times of gap diffusion.}
In: E.~Çinlar, K.~L.~Chung, R.~K.~Getoor (eds), \emph{Seminar on Stochastic Processes, 1981,}
Progress in Probability and Statistics~1, Birkhäuser, Boston, 1981.
\doi{10.1007/978-1-4612-3938-3\textunderscore 3}

\bibitem{kwasnicki-1}
\textsc{Mateusz Kwaśnicki,}
\emph{Harmonic extension technique for non-symmetric operators with completely monotone kernels,}
Calc. Var. Partial. Differ. Equ. 61 (2022), paper no. 202: 1--40.
\doi{10.1007/s00526-022-02308-2}

\bibitem{kwasnicki-2}
\textsc{Mateusz Kwaśnicki,}
\emph{Boundary traces of shift-invariant diffusions in half-plane,}
Ann. Inst. Henri Poincaré Probab. Statist. 59(1) (2023): 411--436.
\doi{10.1214/22-AIHP1250}

\bibitem{kwasnicki-mucha}
\textsc{Mateusz Kwaśnicki, Jacek Mucha,}
\emph{Extension technique for complete Bernstein functions of the Laplace operator,}
J. Evol. Equ. 18(3) (2018): 1341--1379.
\doi{10.1007/s00028-018-0444-4}

\bibitem{molchanov-ostrovskii}
\textsc{Stanislav A.~Molchanov, Evgenii Ostrovskii,}
\emph{Symmetric stable processes as traces of degenerate diffusion processes,}
Theor. Prob. Appl. 14(1) (1969): 128--131.
\doi{10.1137/1114012}

\bibitem{muckenhoupt-stein}
\textsc{Benjamin Muckenhoupt, Elias M.~Stein,}
\emph{Classical expansions and their relation to conjugate harmonic functions,}
Trans. Amer. Math. Soc. 118 (1965): 17--92.
\doi{10.1090/S0002-9947-1965-0199636-9}

\bibitem{schilling-song-vondracek}
\textsc{René L.~Schilling, Renming Song, Zoran Vondraček,}
\emph{Bernstein Functions. Theory and Applications.}
De Gruyter Studies in Mathematics 37, De Gruyter, 2012.
\doi{10.1515/9783110269338}

\bibitem{spitzer}
\textsc{Frank Spitzer,}
\emph{Some theorems concerning 2-dimensional Brownian motion,}
Trans. Amer. Math. Soc. 87 (1958): 187--197.
\doi{10.2307/1993096}

\end{thebibliography}
\end{document}